\documentclass[twoside]{amsart}
\copyrightinfo{}{}

\usepackage{amssymb}
\usepackage{amsmath}
\usepackage{cite}
\usepackage{mathrsfs}

\makeatletter

\address{\bigskip\hfil\begin{tabular}{l@{}}
            School of Mathematics and Statistics F07\\  
            University of Sydney,
            Sydney N.S.W. 2006.\hfill\qquad
                        {\tt mathas@maths.usyd.edu.au}\\  
            Australia.\hfill\qquad
                    {\tt www.maths.usyd.edu.au/u/mathas/}
          \end{tabular}}

\let\atop\@@atop

\def\And{\text{\ and\ }}

\def\For{\text{\ for\ }}

\def\ForSome{\text{\ for some\ }}
\def\If{\text{\ if\ }}
\def\OnlyIf{\text{\ only if\ }}

\def\Otherwise{\text{\ otherwise}}
\def\Unless{\text{\ unless}}
\def\Whenever{\text{\ whenever\ }}
\def\Where{\text{\ where\ }}

\def\Set[#1]#2|#3|{\Big\{\ #2\ \Big| \
            \vcenter{\hsize #1mm\centering#3}\Big\}}

{\catcode`\|=\active
  \gdef\set#1{\mathinner{\lbrace\,{\mathcode`\|"8000%
                                   \let|\midvert #1}\,\rbrace}}
}
\def\midvert{\egroup\mid\bgroup}

\def\Number#1{\refstepcounter{equation}
              \leqno(\theequation)\if*#1%
              \else\def\@currentlabel{{\rm\theequation}}\label{#1}%
              \fi}
\def\Dag{\ifmmode\leqno(\dag)\else$(\dag\)$\fi}
\def\DDag{\ifmmode\leqno(\ddag)\else$(\ddag\)$\fi}

\numberwithin{equation}{section}

\swapnumbers
\newtheorem{Theorem}[equation]{Theorem}
\newtheorem{Proposition}[equation]{Proposition}
\newtheorem{Lemma}[equation]{Lemma}
\newtheorem{Corollary}[equation]{Corollary}

\theoremstyle{remark}
\newtheorem{Assumption}[equation]{Standing assumption}
\newtheorem{Definition}[equation]{Definition}
\newtheorem{Remark}[equation]{Remark}

\newenvironment{Point}[2]%
  {\ifx*#2\let\pointlabel\relax\else\def\pointlabel{#2}\fi
   \refstepcounter{equation}\trivlist
   \item[\hskip\labelsep\bf\theequation
         \ifx\pointlabel\relax\else\space\pointlabel\space\fi]
   \ignorespaces #1
  }{\relax}

\def\Prod{\displaystyle\prod}

\def\){\big)}
\def\({\big(}
\let\>\rangle
\let\<\langle

\let\iso\cong

\let\ss\subseteq
\let\0\varnothing              

\let\realb@r\bar
\let\bar\overline

\let\gedom\trianglerighteq
\let\gdom\vartriangleright

\def\End{\mathop{\rm End}\nolimits}
\def\Hom{\mathop{\rm Hom}\nolimits}

\def\Z{{\mathbb Z}}
\def\Q{{\mathbb Q}}

\let\To\longrightarrow
\def\map#1#2{\,{:}\,#1\!\longrightarrow\!#2}
\def\mapsto{\!\longmapsto\!}

\usepackage{times}

\author{Andrew Mathas} 
\title{Tilting modules for cyclotomic Schur algebras} 

\def\a{\mathfrak a}
\def\b{\mathfrak b}
\def\c{\mathfrak c}
\def\d{\mathfrak d}
\def\s{\mathfrak s}
\def\t{\mathfrak t}
\def\u{\mathfrak u}
\def\v{\mathfrak v}

\def\A{\mathsf A}
\def\B{\mathsf B}
\def\S{\mathsf S}
\def\T{\mathsf T}
\def\U{\mathsf U}
\def\V{\mathsf V}
\def\X{\mathsf X}
\def\Y{\mathsf Y}

\def\O{\mathcal O}
\def\ZZ{\mathcal Z}
\def\HZ{\H_\ZZ}

\def\K{\mathcal K}
\def\HK{\H_\K}

\let\dual\circledast

\def\lamp{{\lambda'}}
\def\H{\mathscr H}
\def\Schur{\mathscr S}
\def\Slam{{\mathcal S}^\lambda}

\let\len\ell
\let\phi\varphi

\def\phiST{\phi_{\S\T}}
\def\phio{\phi_\omega}
\def\phioo{\phi_{\T^\omega\T^\omega}}

\def\Hlambar{\H(\lambda)}
\def\Hlambarp{\H'(\lambda)}

\def\L{L}
\def\M{\mathcal M}
\def\P{P}
\def\W{\Delta}

\def\plam{{\phi_\lambda}}
\def\tlam{{\t^\lambda}}
\def\tllam{{\t_\lambda}}
\def\Tlam{{\T^\lambda}}

\def\tmu{{\t^\mu}}

\def\Sfun{F_\omega}
\def\Mod{\text{--{\bf mod}}}

\def\comp{\operatorname{comp}}
\def\rad{\operatorname{rad}}
\def\res{\operatorname{res}}
\def\Type(#1){\operatorname{Type}(#1)}
\def\Shape(#1){\operatorname{Shape}(#1)}
\def\Sym{\mathfrak S}
\def\rtuple#1{({#1}^{(1)},\dots,{#1}^{(r)})}
\def\Std(#1){{\mathcal T}^{\text{s}}\!(#1)}
\def\SStd(#1,#2){{\mathcal T}^{\text{rs}}_{#2}\!(#1)}
\def\ASStd(#1,#2){{\mathcal T}^{\text{as}}_{#2}\!(#1)}
\def\CStd(#1,#2){{\mathcal T}^{\text{cs}}_{#2}\!(#1)}

\begin{document}
\begin{abstract}
This paper investigates the tilting modules of the
cyclotomic $q$--Schur algebras, the Young modules of the Ariki--Koike
algebras, and the interconnections between them. The main tools used
to understand the tilting modules are contragredient duality, and the
Specht filtrations and dual Specht filtrations of certain permutation
modules.  Surprisingly, Weyl filtrations --- which are in general more
powerful than Specht filtrations --- play only a secondary role.
\end{abstract}
\maketitle

\section{Introduction}

In \cite{DJM:cyc} Dipper, James and the author introduced the
cyclotomic $q$--Schur algebras as another tool for studying the
representations of the Ariki--Koike algebras. These algebras have a
rich and beautiful combinatorial representation theory which closely
resembles that of the $q$--Schur algebras. In particular, the
cyclotomic $q$--Schur algebras are quasi--hereditary and so they have
a theory of tilting modules by Ringel's theorem~\cite{Ringel}.  The
purpose of this paper is to the describe these tilting modules.

A special case of the cyclotomic $q$--Schur algebras are the
$q$--Schur algebras of Dipper and James~\cite{DJ:Schur}.  The
$q$--Schur algebra $\Schur(d,n)$ can be realized as the endomorphism
algebra $\End_{\H(\Sym_n)}(V^{\otimes n})$, where $V$ is the natural
module for the quantum group $U_q(\mathfrak{gl}_d)$ and $\H(\Sym_n)$
is the Iwahori--Hecke algebra of the symmetric group.
Donkin~\cite{Donkin:tilt,Donkin:book} showed that when $d\ge n$ the
tilting modules for $\Schur(d,n)$ are the indecomposable direct
summands of the exterior powers $\wedge^\lambda V
=\wedge^{\lambda_1}V\otimes\dots\otimes\wedge^{\lambda_d}V$, where
$\lambda$ is a partition of $n$.

By definition a cyclotomic $q$--Schur algebra is the $\H$--module
endomorphism algebra of a certain module
$M(\Lambda)=\bigoplus_{\lambda\in\Lambda} M(\lambda)$ of the
Ariki--Koike algebra $\H$. We do not know how to describe $M(\Lambda)$
as a tensor product; however, in essence all of Donkin's results
generalize to the cyclotomic setting --- even though the statements
and proofs do not.  Instead of exterior powers we consider certain
hom--spaces $E(\alpha)=\Hom_\H\(M(\Lambda),N(\alpha)\)$, where
$N(\alpha)$ is something like an induced ``sign representation'' for
the Ariki--Koike algebra. In order to understand these modules, we
work mainly with Specht filtrations and dual Specht filtrations of the
modules $M(\lambda)$ and $N(\alpha)$. Using these filtrations we are
able to show that $M(\lambda)$ and $N(\alpha)$ are both self--dual
$\H$--modules; this implies that $E(\alpha)$ is self--dual. Further,
we can ``lift'' these filtrations to show that $E(\alpha)$ has a Weyl
filtration. Combined, these two results to show that the tilting
modules of the cyclotomic $q$--Schur algebras are the indecomposable
direct summands of the $E(\alpha)$. 

We remark that this description of the tilting modules is valid only
under some mild restrictions on the poset of multipartitions $\Lambda$
and on the defining parameters $Q_1,\dots,Q_r$ for the cyclotomic
$q$--Schur algebra; for the precise statement see Theorem~\ref{tilting modules}. Our first restriction is that $\Lambda$ must contain all
multipartitions of $n$; this is the analogue of the condition $d\ge n$
in Donkin's theorem. The second restriction is that $Q_s\ne0$ for 
any~$s$; this is necessary in order to show that $E(\alpha)$ is
self--dual. The final restriction is that $Q_1,\dots,Q_r$ must all be
distinct; this is needed to force the rank of $E(\alpha)$ to be
independent of the ground ring and the choice of parameters. 

In more detail the contents of the paper are as follows. Section~2
recalls the notation and basic results from the representation theory
of the cyclotomic $q$--Schur algebras and the Ariki--Koike algebras.
Section~3 investigates the Young modules of the Ariki--Koike algebras;
these are the indecomposable direct summands of the modules
$M(\lambda)$ mentioned above. The Young module enjoy all of the
properties of the Young modules of the symmetric groups introduced by
James~\cite{James:trivsour}; they are also closely related to the
tilting modules of the cyclotomic $q$--Schur algebras.  The fourth
section of the paper deals with a duality operation on the category of
$\H$--modules. In the case of the symmetric groups this duality
corresponds to tensoring with the sign representation; applying this
duality to $M(\lambda)$ produces the module $N(\lambda)$. Section~5
studies contragredient duality for the Ariki--Koike algebra; the main
result here is that $M(\lambda)$ and $N(\lambda)$ are both self--dual
$\H$--modules. Finally, building on the previous results, section~6
classifies the tilting modules as the indecomposable direct summands
of the $E(\alpha)$ and section~7 describes the Ringel duals of the
cyclotomic $q$--Schur algebras.

\section{Cyclotomic Schur algebras}

This section is a summary of the definitions and results that we will
need from the representation theory of the cyclotomic $q$--Schur
algebras and the Ariki--Koike algebras. The reader is referred to
\cite{DJM:cyc} for more details.

Fix positive integers $r$ and $n$ and let $\Sym_n$ be the symmetric
group of degree~$n$.  Let $R$ be a commutative ring with $1$ and let
$q,Q_1,\dots,Q_r$ be elements of $R$ such that~$q$ is invertible. The
{\sf Ariki--Koike algebra} $\H=\H_{r,n}$ is the associative unital
$R$--algebra with generators $T_0,T_1,\dots,T_{n-1}$ and relations
$$\begin{array}{r@{\ }l@{\ }ll} (T_0-Q_1)\dots(T_0-Q_r) &=&0, \\
(T_i-q)(T_i+q^{-1}) &=&0,&\text{for $1\le i\le n-1$,}\\
T_0T_1T_0T_1&=&T_1T_0T_1T_0,\\
T_{i+1}T_iT_{i+1}&=&T_iT_{i+1}T_i,&\text{for $1\le i\le n-2$,}\\
T_iT_j&=&T_jT_i,&\text{for $0\le i<j-1\le n-2$.} \end{array}$$ 

The second relation is often written as $(T_i-v)(T_i+1)=0$, for 
$1\le i<n$. This presentation may be turned into the one above by
renormalizing $T_i$ as $v^{-\frac12}T_i$ and setting $q=v^2$. To do
this it is necessary that $v$ have a square root in $R$; however,
every field is a splitting field for $\H$ (because $\H$ is cellular),
so we can adjoin a square root of $v$ without changing the
representation theory of $\H$. We use the presentation above because
it renormalizes the natural inner product on $\H$ and so makes many
formulas nicer later on. We convert the formulas that we need from the
literature without mention.

For $i=1,\dots,n-1$ let $s_i$ be the transposition $(i,i+1)$ in
$\Sym_n$; then $\{s_1,\dots,s_{n-1}\}$ generate $\Sym_n$. If
$w\in\Sym_n$ then $w=s_{i_1}\dots s_{i_k}$ for some $i_j$; if $k$ is
minimal then we say that this expression for $w$ is {\sf reduced} and that
$w$ has {\sf length} $\len(w)=k$. In this case we set $T_w=T_{i_1}\dots
T_{i_k}$; then $T_w$ is independent of the choice of reduced
expression. We also let $L_k=T_{k-1}\dots T_1T_1T_1\dots T_{k-1}$ for
$k=1,2,\dots,n$.  These elements give a basis of~$\H$.

\begin{Point}
{\it}{(Ariki--Koike~\cite[Theorem~3.10]{AK1})} 
The Ariki--Koike algebra $\H$ is free as an $R$--module with basis 
$\set{L_1^{a_1}\dots L_n^{a_n}T_w|w\in\Sym_n
             \And 0\le a_i<r\For 1\le i\le n}$.
\label{AK basis}\end{Point}

Recall that a {\sf composition} of $n$ is sequence
$\sigma=(\sigma_1,\sigma_2,\dots)$ of non--negative integers such that
$|\sigma|=\sum_i\sigma_i=n$; $\sigma$ is a {\sf partition} if in addition
$\sigma_1\ge\sigma_2\ge\cdots$. If $\sigma_i=0$ for all $i>k$ then we
write $\sigma=(\sigma_1,\dots,\sigma_k)$. 

A {\sf multicomposition} of $n$ is an $r$--tuple $\lambda=\rtuple\lambda$
of compositions such that $|\lambda^{(1)}|+\dots+|\lambda^{(r)}|=n$. A
multicomposition $\lambda$ is a {\sf multipartition} if each
$\lambda^{(i)}$ is a partition. If $\lambda$ is a multipartition of
$n$ then we write $\lambda\vdash n$. The diagram $[\lambda]$ of the
multicomposition~$\lambda$ is the set
$[\lambda]=\set{(i,j,s)|1\le\lambda_j^{(s)}\le i\And 1\le s\le r}$.

The set of multicompositions of $n$ is partially ordered by dominance;
that is, if $\lambda$ and $\mu$ are two multicompositions then 
$\lambda$ {\sf dominates} $\mu$, and we write $\lambda\gedom\mu$,~if
$$\sum_{c=1}^{s-1}|\lambda^{(c)}|+\sum_{j=1}^i\lambda^{(s)}_j
   \ge\sum_{c=1}^{s-1}|\mu^{(c)}|+\sum_{j=1}^i\mu^{(s)}_j$$
for $1\le s\le r$ and for all $i\ge 1$. If $\lambda\gedom\mu$ and
$\lambda\ne\mu$ then we write $\lambda\gdom\mu$.

If $\lambda$ is a multicomposition let
$\Sym_\lambda=\Sym_{\lambda^{(1)}}\times\dots\times\Sym_{\lambda^{(r)}}$
be the corresponding Young subgroup of $\Sym_n$. Set
$$x_\lambda=\sum_{w\in\Sym_\lambda}q^{\len(w)}T_w\quad\And\quad
  u_\lambda^+=\prod_{s=2}^r\prod_{k=1}^{a_s}(L_k-Q_s),$$ 
where $a_s=|\lambda^{(1)}|+\dots+|\lambda^{(s-1)}|$ for $2\le s\le r$.
Set $m_\lambda=x_\lambda u_\lambda^+=u_\lambda^+ x_\lambda$ and define
$M(\lambda)$ to be the right ideal $M(\lambda)=m_\lambda\H$ of $\H$.

If $\lambda=\rtuple\lambda$ is a multipartition then a {\sf standard
$\lambda$--tableau} is an $r$--tuple $\t=\rtuple\t$ of standard
tableau which, collectively, contain the integers $1,2,\dots,n$ and
such that $\t^{(c)}$ is a standard $\lambda^{(c)}$--tableau, for 
$1\le c\le r$. Let $\Std(\lambda)$ be the set of standard
$\lambda$--tableau. 

Let $\tlam$ be the standard $\lambda$--tableau with the numbers
$1,2,\dots,n$ entered in order from left to right along its rows. If
$\t$ is any standard $\lambda$--tableau let $d(\t)\in\Sym_n$ be the
unique permutation such that $\t=\tlam d(\t)$. Finally, let
$*\map\H\H$ be the anti--isomorphism given by $T_i^*=T_i^{\phantom*}$,
for $i=0,1,\dots,n-1$, and  set $m_{\s\t}=T_{d(\s)}^*m_\lambda
T_{d(\t)}^{\phantom*}$. 

\begin{Point}
{\it}{(Dipper--James--Mathas \cite[Theorem 3.26]{DJM:cyc})}  
The Ariki--Koike algebra $\H$ is free as an $R$--module with
$($cellular$)$ basis
$\set{m_{\s\t}|\s,\t\in\Std(\lambda)\ForSome
                            \lambda\vdash n}$.
\label{standard basis}\end{Point}

Here, and below, whenever we write expressions involving a pair of
tableaux (such as~$m_{\s\t}$ or $\phiST$), we implicitly assume that
the two tableaux are of the same shape.

The basis $\{m_{\s\t}\}$ is the {\sf standard basis} of $\H$. For
each multipartition $\lambda$ let~$\Hlambar$ be the $R$--submodule
of $\H$ with basis
$\set{m_{\u\v}|\u,\v\in\Std(\mu)\ForSome\mu\gdom\lambda}$; 
then~$\Hlambar$ is a two--sided ideal of $\H$.

Let $S(\lambda)$ be the {\sf Specht module} (or cell module) corresponding
to the multipartition~$\lambda$; that is,
$S(\lambda)\cong(m_\lambda+\Hlambar)\H$, a submodule of $\H/\Hlambar$.
For each $\t\in\Std(\lambda)$ let $m_\t=m_{\tlam\t}+\Hlambar$; then
$S(\lambda)$ is free as an $R$--module with basis
$\set{m_\t|\t\in\Std(\lambda)}$. Further, there is an associative
symmetric bilinear form on $S(\lambda)$ which is determined by 
$$\<m_\s,m_\t\>m_\lambda\equiv m_{\tlam\s}m_{\t\tlam}\mod\Hlambar$$ 
for all $s,\t\in\Std(\lambda)$.  The radical $\rad S(\lambda)$ of this
form is again an $\H$--module, so 
$D(\lambda)=S(\lambda)/\rad S(\lambda)$ is an $\H$--module. When $R$
is a field, $D(\lambda)$ is either $0$ or absolutely irreducible and
all simple $\H$--modules arise uniquely in this way.

We can now give the definition of the cyclotomic $q$--Schur algebras.
A set $\Lambda$ of multicompositions of $n$ is {\sf saturated} if
$\Lambda$ is finite and whenever $\lambda$ is a multipartition such
that $\lambda\gedom\mu$ for some $\mu\in\Lambda$ then
$\lambda\in\Lambda$.  If $\Lambda$ is a saturated set of
multicompositions let $\Lambda^+$ be the set of multipartitions in
$\Lambda$.

\begin{Definition}
Suppose that $\Lambda$ is a saturated set of multipartitions of
$n$.  The {\sf cyclotomic $q$--Schur algebra} with weight poset $\Lambda$
is the endomorphism algebra
$$\Schur(\Lambda)=\End_\H\(M(\Lambda)\),\qquad
     \Where M(\Lambda)=\bigoplus_{\lambda\in\Lambda}M(\lambda).$$
\label{cyclo}\end{Definition}

As we now describe, $\Schur(\Lambda)$ has a basis indexed by pairs of
semistandard tableau. 

A {\sf $\lambda$--tableau of type $\mu$} is a map
$\T\map{[\lambda]}\set{(i,s)|i\ge1\And 1\le s\le r}$ such that
$\mu^{(s)}_i=\#\set{x\in[\lambda]|\T(x)=(i,s)}$ for all $i\ge1$ and
$1\le s\le r$. We think of a $\T$ as being an $r$--tuple
$\T=\rtuple\T$, where $\T^{(s)}$ is the $\lambda^{(t)}$--tableau with
$\T^{(s)}(i,j)=\T(i,j,s)$ for all $(i,j,s)\in[\lambda]$. In this way
we identify the standard tableaux above with the tableaux of type
$\omega=\((0),\dots,(0),(1^n)\)$. If $\T$ is a tableau of type $\mu$
then we write $\Type(\T)=\mu$.

Given two pairs $(i,s)$ and $(j,t)$ write $(i,s)\preceq(j,t)$ if
either $s<t$, or $s=t$ and $i\le j$. 

\begin{Definition}
\label{semistandard} A tableau $\T$ is {\sf (row) semistandard} if,
for $1\le t\le r$, the entries in $\T^{(t)}$ are
\begin{enumerate}
\item weakly increasing along the rows (with respect to $\preceq$);
\item  strictly increasing down columns; and,
\item $(i,s)$ appears in $\T^{(t)}$ only if $s\ge t$. 
\end{enumerate}\noindent 
Let $\SStd(\lambda,\mu)$ be the set of semistandard
$\lambda$--tableau of type $\mu$ and let
$\SStd(\lambda,\Lambda)=\bigcup_{\mu\in\Lambda}\SStd(\lambda,\mu)$
and
$\SStd(\Lambda^+,\mu)=\bigcup_{\lambda\in\Lambda^+}\SStd(\lambda,\mu)$.
\end{Definition}

Usually, we will refer to row semistandard tableaux simply as
semistandard tableaux. Later we will meet column semistandard tableaux
(these are the conjugates of row semistandard tableaux). 

Notice that if $\SStd(\lambda,\mu)$ is non--empty then
$\lambda\gedom\mu$. This observation will be used many times below.

Suppose that $\t$ is a standard $\lambda$--tableau and let $\mu$ be a
multicomposition. Let $\mu(\t)$ be the tableau obtained from $\t$ by
replacing each entry $j$ with $(i,k)$ if $j$ appears in row~$i$ of
$\t^\mu$. The tableau $\mu(\t)$ is a $\lambda$--tableau of type $\mu$;
it is not necessarily semistandard.

If $\S$ and $\T$ are semistandard $\lambda$--tableaux of type $\mu$ 
and $\nu$, respectively, and if~$\t$
is a standard $\lambda$--tableau let 
$$m_{\S\t}=\sum_{\substack{\s\in\Std(\lambda)\\\mu(\s)=\S}}
               q^{\len(d(\s))}m_{\s\t}
\quad\And\quad
m_{\S\T}=\sum_{\substack{\s,\t\in\Std(\lambda)\\\mu(\s)=\S, \nu(\t)=\T}}
           q^{\len(d(\s))+\len(d(\t))}m_{\s\t}.$$
Then we have the following two results.

\begin{Point}
{\it}{\cite[Theorem 4.14]{DJM:cyc}} Suppose that $\mu$ is a
multicomposition of $n$. Then $M(\mu)$ is free as an $R$--module 
with basis
$\set{m_{\S\t}|\S\in\SStd(\lambda,\mu),\t\in\Std(\lambda)
                 \ForSome\lambda\vdash n}$.
\label{sstd basis thm}\end{Point}

For $\S$ and $\T$ as above define $\phiST\in\Schur(\Lambda)$ by
$\phiST(m_\alpha h)=\delta_{\alpha\nu}m_{\S\T}h$, for all $h\in\H$ and
all $\alpha\in\Lambda$. (Here $\delta_{\alpha\nu}$ is the Kronecker
delta; so, $\delta_{\alpha\nu}=1$ if $\alpha=\nu$ and it is zero
otherwise.) Then $\phiST$ belongs to $\Schur(\Lambda)$; moreover,
these elements give us a basis of~$\Schur(\Lambda)$.

\begin{Point}
{\it}{\cite[Theorem~6.6]{DJM:cyc}}
The cyclotomic $q$--Schur algebra $\Schur(\Lambda)$ is free as an
$R$--module with cellular basis
$\set{\phi_{\S\T}|\S,\T\in\SStd(\lambda,\Lambda)
                  \ForSome\lambda\in\Lambda^+}.$
\label{S-basis}\end{Point}

The basis $\{\phiST\}$ is called the {\sf semistandard basis} of
$\Schur(\Lambda)$. Because this basis is cellular the map
$*\map{\Schur(\Lambda)}\Schur(\Lambda)$ which is determined by
$\phi_{\S\T}^*=\phi_{\T\S}$ is an anti--isomorphism of
$\Schur(\Lambda)$. This involution is closely related to the $*$
involution on $\H$; explicitly, if $\phi\map{M(\nu)}M(\mu)$ is an
$\H$--module homomorphism then $\phi^*\map{M(\mu)}M(\nu)$ is the
homomorphism given by $\phi^*(m_\mu h)=\(\phi(\nu)\)^*h$, for all
$h\in\H$.

In order to understand how $\Schur(\Lambda)$ acts on its
representations we need to explain how the multiplication in
$\Schur(\Lambda)$ is determined by the multiplication in $\H$. Suppose
that $\S, \T, \U$ and $\V$ are semistandard tableaux with
$\mu=\Type(\S)$, $\alpha=\Type(\U)$ and $\nu=\Type(\V)$. Then 
$m_{\U\V}=m_\alpha h^\alpha_{\U\V}$, for some
$h^\alpha_{\U\V}\in\H$, and there exist
scalars $r_{\X\Y}\in R$ such that 
$$m_{\S\T}h^\alpha_{\U\V}
  =\sum_{\substack{\X\in\SStd(\Lambda^+,\mu)\\\Y\in\SStd(\Lambda^+,\nu)}}
       r_{\X\Y}m_{\X\Y}                  \Number{S mult}$$
by \cite[Cor~5.17]{DJM:cyc}. Now,
$\phiST\phi_{\U\V}(m_\nu h)=\phiST(m_{\U\V})h
                         =\phiST(m_\alpha)h^\alpha_{\U\V}h
                         =m_{\S\T}h^\alpha_{\U\V}h,$ 
for all $h\in\H$; so (\ref{S mult}) determines the product
$\phiST\phi_{\U\V}$ in~$\Schur(\Lambda)$. Explicitly, we have
$$\phi_{\S\T}\phi_{\U\V}
  =\sum_{\substack{\X\in\SStd(\Lambda^+,\mu)\\\Y\in\SStd(\Lambda^+,\nu)}}
       r_{\X\Y}\phi_{\X\Y},
\Number{Actual S mult}$$
where the $r_{\X\Y}$ are given by (\ref{S mult}).
Note that $\phiST\phi_{\U\V}=0$ if $\Type(\T)\ne\Type(\U)$. In
addition, because $\{\phi_{\S\T}\}$ is a cellular basis, if
$r_{\X\Y}\ne0$ then $\Shape(\X)\gedom\Shape(\S)$, 
with equality only if $\X=\S$. Moreover, if $\X=\S$ then
$r_\Y=r_{\S\Y}$ depends only on $\T$, $\U$ and $\V$; in
particular,~$r_{\S\Y}$ does not depend on $\S$. These details can be
found in \cite{DJM:cyc}; for a complete treatment of the theory of
cellular algebras see \cite{GL,M:ULect}. 

For each multipartition $\lambda\in\Lambda^+$ there is a right
$\Schur(\Lambda)$--module $\W(\lambda)$, called a Weyl module. The Weyl
module $\W(\lambda)$ is the submodule of
$\Hom_\H\(M(\Lambda),S(\lambda)\)$ with basis the set of maps
$\set{\phi_\T|\T\in\SStd(\lambda,\mu),\mu\in\Lambda}$, where
$\phi_\T(m_\alpha h)=\delta_{\alpha\mu}\sum_\t m_\t h$ and the sum is
over those standard $\lambda$--tableaux $\t$ such that $\mu(\t)=\T$.
If $\T$ is a semistandard $\lambda$--tableau and $\phi_{\U\V}$ is
a semistandard basis element then the action of $\Schur(\Lambda)$ on
$\W(\lambda)$ is determined by
$$\phi_\T\phi_{\U\V}
         =\sum_{\Y\in\SStd(\Lambda^+,\nu)} r_{\Y}\phi_\Y,
\Number{W mult}$$
where $r_\Y=r_{\S\Y}$ is determined by~(\ref{S mult}) and
$\nu=\Type(\V)$. (As remarked above, $r_\Y$ is independent of $\S$.)

As with the Specht modules there is an inner product on $\W(\lambda)$
which is determined by
$$\<\phi_\S,\phi_\T\>\plam\equiv\phi_{\Tlam\S}\phi_{\T\Tlam}\mod\Slam,$$
where $\Slam$ is the $R$--submodule of $\Schur(\Lambda)$ with basis
the set of maps $\phi_{\U\V}$ where $\U$ and $\V$ are semistandard
$\mu$--tableaux with $\mu\gdom\lambda$. The quotient module
$\L(\lambda)=\W(\lambda)/\rad \W(\lambda)$ is absolutely irreducible
and $\set{\L(\lambda)|\lambda\in\Lambda^+}$ is a complete set of
non--isomorphic irreducible $\Schur(\Lambda)$--modules.

Recall that $\omega=\rtuple\omega$ is the multipartition with
$\omega^{(r)}=(1^n)$ and $\omega^{(s)}=(0)$ for $1\le s<r$. From the
definitions, $m_\omega=1$ and $\phio=\phioo$ is the identity map 
on~$\H$; so, $\H=M(\omega)$. In particular, $\phio$ is an idempotent in
$\Schur(\Lambda)$ and it is easy to see that
$\H\cong\phio\Schur(\Lambda)\phio$ whenever $\omega\in\Lambda$.  

For an algebra $A$ let $A\Mod$ be the category of finite dimensional
right $A$--modules. As noted in \cite{JM:cyc-Schaper}, standard
arguments show that there is a functor 
$$\Sfun\map{\Schur(\Lambda)\Mod}\H\Mod; M\mapsto M\phio$$
which has the following properties.

\begin{Point}
{\it}{(The cyclotomic Schur functor~\cite{JM:cyc-Schaper})} 
Suppose that $R$ is a field and that $\omega\in\Lambda$. Let
$\lambda\in\Lambda^+$. Then, as right $\H$--modules,
\begin{enumerate}
\item $\Sfun(\W(\lambda))\iso S(\lambda)$;
\item $\Sfun(\L(\lambda))\iso D(\lambda)$.
\end{enumerate} \noindent 
Furthermore, if $D(\mu)\ne0$ then 
$[\W(\lambda):\L(\mu)]=[S(\lambda):D(\mu)]$.
\label{Schur functor}\end{Point}


\section{Cyclotomic Young modules}

For each multicomposition $\mu$ of $n$ let $\phi_\mu$ be the identity
map on $M(\mu)$.  This section describes the indecomposable summands
of $M(\mu)$. We approach this question by considering the right
$\Schur(\Lambda)$--modules $\M(\mu)=\Hom_\H\(M(\Lambda),M(\mu)\)$.

All of these results in this section about the modules $\M(\mu)$ apply
without restriction on~$\Lambda$; however, whenever we apply the Schur
functor we implicitly assume that $\omega\in\Lambda$.

\begin{Proposition}
Suppose that $\mu\in\Lambda$. Then the following hold.
\begin{enumerate}
\item $\M(\mu)$ is free as an
$R$--module with basis
$$\set{\phi_{\S\T}|\S\in\SStd(\lambda,\mu),\T\in\SStd(\lambda,\nu)
\ForSome\nu\in\Lambda\And\lambda\in\Lambda^+}.$$ 
\item $\M(\mu)\cong\phi_\mu\Schur(\Lambda)$ as right
$\Schur(\Lambda)$--modules; in particular, $\M(\mu)$ is projective.
\item As $\H$--modules, $M(\mu)\cong\Sfun(\M(\mu))$.  
\end{enumerate}\label{M^mu S-basis}
\end{Proposition}

\begin{proof}
Part (i) is just a restatement of (\ref{S-basis}). For
(ii), note that if $\S\in\SStd(\lambda,\alpha)$, for some
$\lambda\in\Lambda^+$ and $\alpha\in\Lambda$, then
$\phi_\mu\phiST=\delta_{\alpha\mu}\phiST$ for all
$\T\in\SStd(\lambda,\Lambda)$. Hence,
$\phi_\mu\Schur(\Lambda)=\M(\mu)$ by part~(i). As $\phi_\mu$ is an
idempotent this also shows that $\M(\mu)$ is a projective
$\Schur(\Lambda)$--module. Finally, by part (i)
again, the $\H$--module $\Sfun(\M(\mu))=\Hom_\H\(\H,M(\mu)\)$ is free
as an $R$--module with basis
$\set{\phi_{\S\t}|\S\in\SStd(\lambda,\mu),\t\in\Std(\lambda)
                        \ForSome\lambda\in\Lambda^+}$.
Hence, by (\ref{sstd basis thm}), $\Sfun(\M(\mu))\cong M(\mu)$, where the
isomorphism is given by the $R$--linear map determined by
$\phi_{\S\t}\mapsto m_{\S\t}=\phi_{\S\t}(m_\mu)$ for all
$\S\in\SStd(\lambda,\mu)$ and $\t\in\Std(\lambda)$.
\end{proof}

An $\Schur(\Lambda)$--module $X$ has a {\sf Weyl filtration} if it has an
$\Schur(\Lambda)$--module filtration
$$X=X_1\supset\dots\supset X_k\supset X_{k+1}=0$$
such that $X_i/X_{i+1}\cong\W(\lambda_i)$ for some multipartition
$\lambda_i\in\Lambda^+$, for $1\le i\le k$. Since each Weyl module
$\W(\lambda)$ has simple head $\L(\lambda)$ and
$[\rad\W(\lambda):L(\mu)]\ne0$ only if $\lambda\gedom\mu$ the equivalence
classes of the Weyl modules are a basis of the Grothendieck group of
$\Schur(\Lambda)$; consequently, when $R$ is a field the filtration
multiplicities $$[X:\W(\lambda)]=\#\set{1\le i\le
k|X_i/X_{i+1}\cong\W(\lambda)}$$ are independent of the choice of
filtration. Finally, note that if $X$ has a Weyl filtration as above
then
$\Sfun(X)=\Sfun(X_1\supset\dots\supset\Sfun(X_k)\supset\Sfun(X_{k+1})=0$
is a {\sf Specht filtration} of~$\Sfun(X)$ by (\ref{Schur functor}); that
is, 
$\Sfun(X_i)/\Sfun(X_{i+1})\cong\Sfun(X_i/X_{i+1})\cong S(\lambda_i)$,
for $1\le i\le k$.

\begin{Lemma}
Suppose that $\mu\in\Lambda$. Then $\M(\mu)$ has a Weyl
filtration 
$$\M(\mu)=\M_1\supset \M_2\supset \dots\supset \M_k\supset\M_{k+1}=0$$
and there exist multipartitions $\lambda_1,\dots,\lambda_k$ such that
$\M_i/\M_{i+1}\cong\W(\lambda_i)$, for $i=1,\dots,k$.  Moreover,
if $\lambda_i\gdom\lambda_j$ then $i>j$ and 
$\#\set{1\le i\le k|\lambda_i=\lambda}=\#\SStd(\lambda,\mu)$ for each
multipartition $\lambda$. Hence,
$[\M(\mu){:}\W(\lambda)]=\#\SStd(\lambda,\mu)$ when $R$
is a field.
\label{Weyl filtration}
\end{Lemma}

\begin{proof}
By Proposition~\ref{M^mu S-basis},
$\set{\phi_{\S\T}|\S\in\SStd(\lambda,\mu),\T\in\SStd(\lambda,\nu)
                        \ForSome\nu\in\Lambda\And\lambda\in\Lambda^+}$
is a basis of~$\M(\mu)$. Let $\{\S_1,\dots,\S_k\}$ be the set of
semistandard tableaux of type $\mu$ ordered so that
$i>j$ whenever $\Shape(\S_i)\gdom\Shape(\S_j)$. Let
$\lambda_i=\Shape(\S_i)$, for $1\le i\le k$. Notice that
$\lambda_i\gedom\mu$, for all $i$, since $\SStd(\lambda_i,\mu)\ne\0$.

Fix an integer $i$, with $1\le i\le k$, and let $\M_i$ be the 
$R$--submodule of $\M(\mu)$ with basis
$$\set{\phi_{\S_j\T}|\lambda_j\gedom\lambda_i\And
                     \T\in\SStd(\lambda_j,\Lambda)}.$$
Then $\M_i$ is an $\Schur(\Lambda)$--module by the remarks after
(\ref{Actual S mult}). Further, there is an isomorphism of
$\Schur(\Lambda)$--modules $\W(\lambda_i)\cong \M_i/\M_{i+1}$ given by
$\phi_\T\mapsto\phi_{\S_i\T}+\M_{i+1}$, for
$\T\in\SStd(\lambda_i,\Lambda)$, because 
$\mathcal S^{\lambda_i}\cap \M_i\subseteq \M_{i+1}$. 
\end{proof}

The simple $\Schur(\Lambda)$--modules $\L(\lambda)$ are indexed by the
multipartitions $\lambda\in\Lambda^+$. Fix a set $\set{
\P(\lambda)|\lambda\in\Lambda^+}$ of principal indecomposable
$\Schur(\Lambda)$--modules where $\P(\lambda)$ is the projective cover
of $\L(\lambda)$.

\begin{Proposition}
Suppose that $R$ is a field and let $\mu$ be a
multipartition of $n$.  Then 
$$\M(\mu)
   \cong \P(\mu)\oplus\bigoplus_{\lambda\gdom\mu}c_{\lambda\mu} \P(\lambda)$$
for some non--negative integers $c_{\lambda\mu}$.
\label{M^mu S-filtration}
\end{Proposition}

\begin{proof}
By Proposition~\ref{M^mu S-basis}, $\M(\mu)$  is a projective
$\Schur(\Lambda)$--module; therefore, there exist non--negative
integers $c_{\lambda\mu}$ such that 
$\M(\mu)\cong\bigoplus_\lambda c_{\lambda\mu} \P(\lambda)$. Now, by
\cite[Theorem~3.7]{GL} (or, more explicitly,
\cite[Lemma~2.19]{M:ULect}), each $\P(\lambda)$ has a Weyl filtration
in which $\W(\lambda)$ appears with multiplicity~$1$. On the other
hand, Lemma~\ref{Weyl filtration} $\M(\mu)$ has a Weyl filtration in which
the Weyl module $\W(\lambda)$ is a subquotient only if
$\SStd(\lambda,\mu)$ is non--empty; that is, if~$\lambda\gedom\mu$.
Hence, $c_{\lambda\mu}\ne0$ only if $\lambda\gedom\mu$. 

It remains to show that $c_{\mu\mu}=1$. First, observe that by
Lemma~\ref{Weyl filtration} $\W(\mu)$ is a top composition factor of
$\M(\mu)$; consequently, $\L(\mu)$ is also top composition factor of
$\M(\mu)$. On the other hand, $\L(\mu)$ is a top composition factor of
$\P(\lambda)$ if and only if $\lambda=\mu$: hence, $\P(\mu)$ is a
direct summand of $\M(\mu)$ and $c_{\mu\mu}\ge1$. Therefore, by
Lemma~\ref{Weyl filtration},
$$1=[\M(\mu):\W(\mu)]
   =\sum_\lambda c_{\lambda\mu}[\P(\lambda):\W(\mu)]
   \ge c_{\mu\mu}\ge1.$$  
We must have equality throughout; so $c_{\mu\mu}=1$ and the
Proposition follows.
\end{proof}

Suppose that $\S\in\SStd(\lambda,\mu)$ and $\T\in\SStd(\lambda,\nu)$.
Then, since $\M(\mu)=\Hom_\H\(M(\Lambda),M(\mu)\)$, we
can define an $\Schur(\Lambda)$--module homomorphism
$\Phi_{\S\T}\map{\M(\nu)}\M(\mu)$ by
$\Phi_{\S\T}(f)=\phi_{\S\T}f$ for all
$f\in\M(\nu)$. In fact, as $\S$ and $\T$ run over 
$\SStd(\lambda,\mu)$ and $\SStd(\lambda,\nu)$, respectively, these 
maps give a basis of $\Hom_\H(\M(\nu),M(\mu)\)$. 

\begin{Lemma}
Suppose that $\mu\in\Lambda$. Then
$\Hom_{\Schur(\Lambda)}(\M(\nu),\M(\mu))$ is free as an $R$--module with
basis 
$\set{\Phi_{\S\T}|\S\in\SStd(\lambda,\mu),\T\in\SStd(\lambda,\nu)
           \ForSome\lambda\vdash n}. $
\label{psi basis}
\end{Lemma}

\begin{proof}
By definition each of the maps $\Phi_{\S\T}$ belongs to
$\Hom_{\Schur(\Lambda)}(\M(\nu),\M(\mu))$ and they are certainly linearly
independent. It remains to check that these homomorphisms span
$\Hom_{\Schur(\Lambda)}(\M(\nu),\M(\mu))$. Now, if
$f\in\Hom_{\Schur(\Lambda)}(\M(\nu),\M(\mu))$ then there exist
$a_{\S\T}\in R$ such that
$f(\phi_\nu)=\sum_{\S\T}a_{\S\T}\phiST$ by Proposition~\ref{M^mu S-basis}. Hence, 
$f=\sum a_{\S\T}\Phi_{\S\T}$ and the Lemma is proved.
\end{proof}

For each multipartition $\lambda$ let $Y(\lambda)=\Sfun(
\P(\lambda))$. If $\mu\in\Lambda$ is a multipartition then
$M(\mu)\cong\Sfun\(\M(\mu)\)$ by Proposition~\ref{M^mu S-basis}; therefore, by
Proposition~\ref{M^mu S-filtration},
$$M(\mu)\cong Y(\mu)\oplus
        \bigoplus_{\lambda\gdom\mu}c_{\lambda\mu}Y(\lambda).
\Number{Young module defn}$$
As remarked above, $P(\lambda)$ has a Weyl filtration. Therefore,
$Y(\lambda)$ has a Specht filtration; in particular, $Y(\lambda)\ne0$.
Following James~\cite{James:trivsour},  we call $Y(\lambda)$ a 
{\sf Young module} of~$\H$.

\begin{Theorem}
\label{Young modules} Suppose that $R$ is a field and let $\mu$ 
be a multipartition of $n$. Then the following hold.
\begin{enumerate}
\item Each $Y(\mu)$ is an indecomposable $\H$--module;
\item If $\lambda$ is another multipartition of $n$ then
$Y(\lambda)\cong Y(\mu)$ if and only if $\lambda=\mu$; and,
\item The Young module $Y(\mu)$ has a Specht filtration
$$Y(\mu)=Y_1\supset\dots\supset Y_k\supset Y_{k+1}=0$$ 
with $Y_i/Y_{i+1}\cong S^{\lambda_i}$, for some multipartitions
$\lambda_1,\dots,\lambda_k$.
\item The number of $\lambda_i$ equal to $\lambda$ is  
the decomposition multiplicity $[\W(\lambda):\L(\mu)]$.
\end{enumerate}
\end{Theorem}

\begin{proof}
First note that (\ref{S-basis}) and Lemma~\ref{psi basis} show that
$$\Hom_{\Schur(\Lambda)}(\M(\mu),\M(\lambda))
          \cong\Hom_\H(M(\lambda),M(\mu))$$
as $R$--modules; explicitly, the isomorphism is given by
$\Phi_{\S\T}\mapsto\Sfun(\Phi_{\S\T})=\phiST$. Therefore, $\Sfun$
induces an injective map
$\Hom_{\Schur(\Lambda)}(\P(\mu),\P(\lambda))
             \hookrightarrow\Hom_\H(Y(\lambda),Y(\mu))$.
Now $Y(\lambda)$ is a direct summand of $M(\mu)$ so any map from
$Y(\lambda)$ to $Y(\mu)$ can be extended to a map in
$\Hom_\H(M(\lambda),M(\mu))$; hence,
$\Hom_{\Schur(\Lambda)}(\P(\mu),\P(\lambda))$ and 
$\Hom_\H(Y(\lambda),Y(\mu))$
are isomorphic $R$--modules. 

In the special case where $\lambda=\nu$ the last paragraph says that
$\End_{\Schur(\Lambda)}(\P(\mu))$ and $\End_\H(Y(\mu))$ are isomorphic
rings.  This proves~(i) as $\End_{\Schur(\Lambda)}(\P(\mu))$ is
a local ring because $P(\mu)$ is indecomposable. Similarly, part (ii)
follows because if $Y(\mu)\cong Y(\lambda)$ then
$\Hom_\H(Y(\mu),Y(\lambda))$ contains an isomorphism and this lifts to
give an isomorphism $\P(\mu)\cong \P(\lambda)$, so $\lambda=\mu$.

We now prove (iii). Recall from the proof of Proposition~\ref{M^mu S-filtration}
that $\P(\mu)$ has a Weyl filtration
$\P(\mu)=P_1\supset\dots\supset P_k\supset P_{k+1}=0$. Moreover, for 
each multipartition $\lambda$,
$$\#\set{1\le i\le k|P_i/P_{i+1}\cong\W(\lambda)}
      =[P(\mu):\W(\lambda)]=[\W(\lambda):\L(\mu)],$$ 
where the last equality follows from \cite[Lemma~2.19]{GL} (the
cellular algebra analogue of the Brauer--Nesbitt $cde$--triangle).
Setting $Y_i=\Sfun(P_i)$, and using (\ref{Schur functor}), gives a
filtration of~$Y(\mu)$ with the required properties. Notice that
$\Sfun(\W(\lambda))\cong S(\lambda)$ for all~$\lambda$; therefore,
even though $\Sfun(\L(\lambda))=0$ when $D(\lambda)=0$ the
multiplicities in the Specht filtration of~$Y(\mu)$ are preserved.
\end{proof}

In part (iii) we can do slightly better because the arguments of
\cite{GL,M:ULect} show that $\P(\mu)$ can be filtered so that each of
the quotients is isomorphic to a direct sum of $[\W(\lambda):\L(\mu)]$
copies of the Weyl module $\W(\lambda)$.

Let $(K,\O,R)$ be a modular system (with parameters). That is,
$\O\subset K$ is a discrete valuation ring with residue field $R$ and
we choose parameters $\hat q,\hat Q_1,\dots,\hat Q_r$ in~$\O$ so that
the Ariki--Koike algebra $\H_K$ over $K$ with parameters $\hat q,\hat
Q_1,\dots,\hat Q_r\in\O$ is semisimple and $\pi(\hat q)=q$ and
$\pi(\hat Q_s)=Q_s$, for $1\le s\le r$, where $\pi\map\O R$ is the
canonical projection map. Let $\H_{\O}$ be the Ariki--Koike algebra
with parameters $\hat q,\hat Q_1,\dots,\hat Q_r\in\O$; then
$\H_K\cong\H_\O\otimes_\O K$ and $\H=\H_R\cong\H_\O\otimes_O R$.

Let $Y$ be an $\H_R$--module with an $\O$ lattice;
that is, an $\O$--free $\H_\O$--module $Y_\O$ such that
$Y\cong Y_\O\otimes_\O R$. Suppose that $Y_\O$ has a Specht filtration
$$Y_\O=Y_{\O,1}\supset\dots\supset Y_{\O,k}\supset0$$
and set $Y_i=Y_{\O,i}\otimes_\O R$, for all $i$. Then 
$Y\supset Y_1\supset\dots\supset Y_k\supset0$ is a Specht filtration of 
$Y$. In this case for any multipartition $\lambda$ we define 
$$[Y:S(\lambda)]=\dim_K\Hom_{\H_K}(Y_\O\otimes_\O K,S(\lambda)_K).$$
Then $[Y:S(\lambda)]$ is independent of the choice of lattice $Y_\O$
and the choice of filtration ($Y$ is a modular reduction of 
$Y_K=Y_\O\otimes_\O K$ and $Y_K$ is independent of these choices being
semisimple). 

As Theorem~\ref{Young modules}(iii) holds for all rings we can rephrase
Theorem~\ref{Young modules}(iv) as follows.

\begin{Corollary}
Suppose that $R$ is a field and let $\lambda$ and $\mu$ be
multipartitions of $n$. Then
$$[Y(\mu):S(\lambda)]=[\W(\lambda):L(\mu)].$$
\label{Young mult}
\end{Corollary}

Note that we cannot just define $[Y(\mu):S(\lambda)]$ to be equal to
the number of subquotients which are isomorphic to $S(\lambda)$ in a
Specht filtration of $Y(\mu)$ because it can happen that
$S(\lambda)\cong S(\nu)$ even though $\lambda\ne\nu$. This is why we
have to introduce a modular system.

If $\mu=\rtuple\mu$ is a multicomposition of $n$ let
$\vec\mu=\rtuple{\vec\mu}$ be the unique multipartition of $n$
such that $\vec\mu^{(i)}$ is the partition obtained from $\mu^{(i)}$
by reordering its parts. The following result is needed in
\cite{DM:Morita}.

\begin{Corollary}
Suppose that $R$ is a field and let $\mu$ be a multicomposition 
of $n$. Then
$$M(\mu)\cong Y(\vec\mu)
   \oplus\bigoplus_{\lambda\gdom\vec\mu}c_{\lambda\vec\mu}Y(\lambda)$$
where the integers $c_{\lambda\vec\mu}$ are as in 
Proposition~\ref{M^mu S-filtration}.
\label{M^mu summands}
\end{Corollary}

\begin{proof}
If $\mu$ is a multipartition then this is just a restatement of
(\ref{Young module defn}), so suppose that $\mu$ is not a
multipartition. Then $\Sym_\mu$ and $\Sym_{\vec\mu}$ are conjugate
subgroups of $\Sym_n$; therefore we can find a permutation
$d\in\Sym_n$ such that $\Sym_\mu=d^{-1}\Sym_{\vec\mu}d$ and $\tmu
d^{-1}$ and $\t^{\vec\mu}d$ are both row standard (see, for example,
\cite[Lemma~3.10]{M:ULect}).  For this $d$ we have $T_d
m_\mu=m_{\vec\mu}T_d$ (by \cite[2.1(iv)]{DJM:cyc}); consequently,
$M(\mu)\cong T_d^{-1}M(\vec\mu)\cong M(\vec\mu)$ as right
$\H$--modules. The general case now follows from 
(\ref{Young module defn}).
\end{proof}


\section{Twisted cyclotomic Schur algebras}

The Ariki--Koike algebra $\H=\H_{r,n}$ has (at most) $2r$ one dimensional
characters; namely, the $R$--linear maps
$\chi_{s,\alpha}\map{\H_{r,n}}R$, for $1\le s\le r$ and 
$\alpha\in\{q,-q^{-1}\}$, which are determined by
$\chi_{s,\alpha}(T_0)=Q_s$ and
$\chi_{s,\alpha}(T_i)=\alpha$, for $1\le i<n$. The character
$\chi_{s,\alpha}$ is afforded by the Specht module 
$S^{\lambda_{s,\alpha}}$ where 
$$\lambda^{(t)}_{s,\alpha}
       =\begin{cases}(n),&\If s=t\And\alpha=q,\\
                     (1^n),&\If s=t\And\alpha=-q^{-1},\\
                     (0),&\Otherwise.
\end{cases}$$
Clearly, $S^{\lambda_{s,\alpha}}\cong S^{\lambda_{t,\beta}}$ if and
only if $(Q_s,\alpha)=(Q_t,\beta)$. When $\H$ is semisimple all of
these representations are pairwise non--isomorphic.

Given any $\H$--module $M$ we can use the character
$\chi_{s,\alpha}$ to twist the $\H$--action to give a new $\H$--module
$M_{s,\alpha}$ on which $h\in\H$ acts as $\chi_{s,\alpha}(h)h$. By
considering characters, in the semisimple case the effect of this
operation on the Specht modules amounts to a cyclic permutation of the
components of the corresponding multipartitions, and taking conjugates
when $\alpha=-q^{-1}$. In contrast, when $\H$ is not semisimple the
twisted Specht module $S^\lambda_{s,\alpha}$ is not necessarily
isomorphic to another Specht module. 

This section investigates what happens when we twist modules by the
`sign representation' $\chi_{r,-q^{-1}}$ of $\H$. These twisted modules
play a key role in understanding the tilting modules of the cyclotomic
Schur algebras.

Let $\ZZ=\Z[\dot q,\dot q^{-1},\dot Q_1,\dots,\dot Q_r]$, where
$\dot q,\dot Q_1,\dots,\dot Q_r$ are indeterminates over~$\Z$, and let
$\HZ$ be the Ariki--Koike algebra over $\ZZ$ with parameters 
$\dot q,\dot Q_1,\dots,\dot Q_r$. The relations of $\H$ imply that
$\H$ has a~$\Z$--algebra involution ${}'$ which is determined by
$$T_i'=T_i,\quad
  \dot q'=-\dot q^{-1},\quad\And\quad
  \dot Q_s'=\dot Q_{r-s+1}, $$
for $0\le i<n$ and $1\le s\le r$. Then $L_k'=L_k$ and $T_w'=T_w$, for
all $1\le k\le n$ and $w\in\Sym_n$. We emphasize that the involution
${}'$ is only defined generically (i.e.~over $\ZZ$) and that $\H$ does
not have a corresponding involution when $R$ is not a free
$\ZZ$--module under specialization. Nevertheless, specialization
arguments will allow us to transport the effects of ${}'$ into $\H_R$.

Suppose that $\lambda$ is a multicomposition and define
$$y_\lambda=\sum_{w\in\Sym_\lambda}(-\dot q)^{-\len(w)}T_w\quad\And\quad
  u_\lambda^-=\prod_{s=1}^{r-1}\prod_{k=1}^{b_s}(L_k-\dot Q_s),$$ 
where
$b_s=|\lambda^{(s+1)}|+\dots+|\lambda^{(r)}|$ for $2\le k\le r$. 
Then $y_\lambda=(x_\lambda)'$ and $u^-_\lambda=(u^+_\lambda)'$; in
particular, it follows that
$y_\lambda u_\lambda^-=u_\lambda^-y_\lambda$. Set
$n_\lambda=y_\lambda u_\lambda^-$ and, if  $\s$ and $\t$ are
standard $\lambda$--tableaux, define
$n_{\s\t} =T_{d(\s)}^* n_\lambda T_{d(\t)}$;
then $n_{\s\t}=(m_{\s\t})'$. Therefore, because ${}'$ is a $\Z$--algebra
involution, $\{n_{\s\t}\}$ is a cellular
basis of $\HZ$ by (\ref{standard basis}).

Returning to the general case, any ring $R$ with a choice of
parameters $q,Q_1,\dots,Q_r$ is naturally a~$\ZZ$--module under
specialization: that is, $\dot q$ acts on $R$ as multiplication 
by~$q$, and~$\dot Q_s$ acts as multiplication by $Q_s$, for 
$1\le s\le r$.  Moreover, because $\H$ is $R$--free this induces a
isomorphism of $R$--algebras $\H_R\cong\HZ\otimes_\ZZ R$ via
$T_i\mapsto T_i\otimes 1_R$, for $0\le i<n$. We say that $\H_R$ is a
{\sf specialization} of~$\HZ$.

Hereafter, we drop the distinction between $q$ and $\dot q$, and 
$\dot Q_s$ and $Q_s$, and we identify the algebras $\H=\H_R$ and
$\HZ\otimes_\ZZ R$ via the isomorphism $T_i\mapsto T_i\otimes 1_R$
above. Thus, we have elements $y_\lambda$, $u_\lambda^-$ and
$n_{\s\t}$ in $\H$ and by (\ref{standard basis}), and the specialization
argument above, we have the following.

\begin{Point}
{\it}{(Du--Rui \cite[2.7]{DuRui:branching})}  
The Ariki--Koike algebra $\H$ is free as an $R$--module with
cellular basis
$\set{n_{\s\t}|\s,\t\in\Std(\lambda)\ForSome\lambda\vdash n}$.
\label{costandard basis}\end{Point}

Since $\{n_{\s\t}\}$ is a cellular basis it gives us a second
collection of cell modules for~$\H$; namely, for each multipartition
$\lambda$ define the {\sf dual Specht module} $S'(\lambda)$ to be the
right $\H$--module $(n_\lambda+\Hlambarp)\H$, where
$\Hlambarp=\(\Hlambar\)'$ is the two--sided ideal of $\H$ with basis
$n_{\u\v}$ with $\Shape(\u)=\Shape(\v)\gdom\lambda$.  Then
$S'(\lambda)$ is $R$--free with basis $\set{n_\t|\t\in\Std(\lambda)}$,
where $n_\t=n_{\tlam\t}+\Hlambarp$. This terminology is justified 
in Corollary~\ref{dual Spechts} below which shows that $S'(\lambda)$ is
isomorphic to the contragredient dual of $S(\lambda')$, where
$\lambda'$ is the multipartition conjugate to~$\lambda$. We remark
when $\H$ is semisimple a straightforward calculation using characters
shows that $S'(\lambda)\cong S(\lambda')_{r,-q^{-1}}$. 

Let $D'(\lambda)=S'(\lambda)/\rad S'(\lambda)$, where $\rad
S'(\lambda)$ is the radical of the bilinear form on $S'(\lambda)$; the
form is defined in terms of the structure constants of the cellular
basis $\{n_{\s\t}\}$. Once again, the theory of cellular algebras says
that the non--zero $D'(\lambda)$ are a complete set of pairwise
non--isomorphic irreducible $\H$--modules.

For any multicomposition $\mu$ let $N(\mu)=n_\mu\H$. If $\S$ is a
semistandard $\lambda$--tableau of type~$\mu$ and $\t$ is a standard
$\lambda$--tableau define
$$n_{\S\t}=\sum_{\substack{\s\in\Std(\lambda)\\\mu(\s)=\S}} 
                   (-q)^{-\len(d(\s))}n_{\s\t}.$$
From the definitions, $n_{\S\t}=m_{\S\t}'$ in $\HZ$; therefore,
(\ref{sstd basis thm}) and the usual specialization argument show that
the following holds.

\begin{Corollary}
Suppose that $\mu$ is a multicomposition of $n$. Then $N(\mu)$ is
free as an $R$--module with basis
$\set{n_{\S\t}|\S\in\SStd(\lambda,\mu),\t\in\Std(\lambda)\ForSome
                            \lambda\vdash n}$.
\label{N sstd basis thm}\end{Corollary}

Just as in \cite[Cor.~4.15]{DJM:cyc}, this implies that $N(\mu)$ has a
dual Specht filtration in which the number of subquotients equal to
$S'(\lambda)$ is $\#\SStd(\lambda,\mu)$. This filtration can also be
obtained by specializing the corresponding Specht filtration of
the $\HZ$--module $M(\mu)$.

Mirroring Definition~\ref{cyclo}, if $\Lambda$ is a saturated set of
multicompositions $\Lambda$ define the {\sf twisted cyclotomic
$q$--Schur algebra} to be the endomorphism algebra
$$\Schur'(\Lambda)=\End_\H\(N(\Lambda)\),\qquad
     \Where N(\Lambda)=\bigoplus_{\mu\in\Lambda}N(\mu).$$

If $\S\in\SStd(\lambda,\mu)$ and $\T\in\SStd(\lambda,\nu)$ are
semistandard tableaux let 
$$n_{\S\T}
  =\sum_{\substack{\s,\t\in\Std(\lambda)\\\mu(\s)=\S,\,\nu(\t)=\T}} 
         (-q)^{-\len(d(\s))-\len(d(\t))}n_{\s\t}.$$ 
Now define the homomorphism $\phiST'\in\Schur'(\Lambda)$ by
$\phiST'(n_\alpha h)=\delta_{\alpha\nu}n_{\S\T}h$, for all $h\in\H$
and all $\alpha\in\Lambda$. Then $\phiST'$ belongs to
$\Schur'(\Lambda)$.

Write $\Schur'(\Lambda)_\ZZ$ for the twisted cyclotomic Schur algebra
over $\ZZ$. Similarly, we write $M(\mu)_\ZZ$, $N(\mu)_\ZZ,\dots$
whenever we have a free $R$--module whose rank is independent of $R$,
$q$ and $Q_1,\dots,Q_r$.

\begin{Proposition}
\label{S' basis}
Suppose that $\Lambda$ is a saturated set of multicompositions.
\begin{enumerate}
\item The twisted cyclotomic $q$--Schur algebra $\Schur'(\Lambda)$ 
is free as an $R$--module with cellular basis
$$\set{\phiST'|\S\in\SStd(\lambda,\mu),\T\in\SStd(\lambda,\nu)
       \ForSome\mu,\nu\in\Lambda\text{\ and some\ }\lambda\in\Lambda^+}. $$
Consequently, $\Schur'(\Lambda)\cong\Schur'(\Lambda)_\ZZ\otimes_\ZZ R$.
\item The twisted cyclotomic Schur algebra $\Schur'(\Lambda)$ is 
quasi--hereditary.
\item The $R$--algebras $\Schur'(\Lambda)$ and $\Schur(\Lambda)$ are
canonically isomorphic.
\end{enumerate}
\end{Proposition}

\begin{proof}
Using Corollary~\ref{N sstd basis thm}, an easy modification of
the argument of \cite[Prop.~6.3]{DJM:cyc} shows that 
$\set{n_{\S\T}|\S\in\SStd(\lambda,\mu), \T\in\SStd(\lambda,\nu)
               \ForSome\lambda\vdash n}$ 
is a basis of $N(\nu)^*\cap N(\mu)$. Part (i) nows follows exactly as
in the proof (\ref{S-basis}); see \cite[Theorem~6.6]{DJM:cyc}. In
particular, notice that because 
$\Schur'(\Lambda)\cong\Schur'(\Lambda)_\ZZ\otimes_\ZZ R$ we can now use specialization arguments. 

Part (ii) follows from (i) using the argument of
\cite[Cor.~6.18]{DJM:cyc}; alternatively, it may be deduced by the
specialization of a hereditary chain of the algebra 
$\Schur(\Lambda)_\ZZ$. 

Finally, (iii) follows because when $R=\ZZ$,
$$\Hom_\H\(M(\nu),M(\mu)\)\cong\Hom_\H\(M(\nu)',M(\mu)'\)
              =\Hom_\H\(N(\nu),N(\mu)\);$$
explicitly, the isomorphism is given by $\phiST\mapsto\phiST'$, for
semistandard tableaux $\S$ and $\T$. As
$\Schur(\Lambda)\cong\Schur(\Lambda)_\ZZ\otimes R$ and 
$\Schur'(\Lambda)\cong\Schur'(\Lambda)_\ZZ\otimes R$ this implies the
general case.
\end{proof}

Let $\W'(\lambda)$ and $L'(\lambda)$, respectively, be the Weyl
modules and simple modules of $\Schur'(\Lambda)$; these are defined in
exactly the same way as the corresponding modules for
$\Schur(\Lambda)$. As in (\ref{Schur functor}), if $\omega\in\Lambda$
then there is a functor 
$$\Schur'(\Lambda)\Mod\To\H\Mod; N\mapsto N\phi_\omega',$$ 
where $\phi_\omega'$ is the identity map on $\H$.  Because
$\phi_\omega'=\phi_\omega$ we abuse notation and again denote this
functor by $\Sfun$. As in (\ref{Schur functor}), we have
$\Sfun(\W'(\lambda))\cong S'(\lambda)$, $\Sfun(L'(\lambda))\cong
D'(\lambda)$ and $[\W'(\lambda):L'(\mu)]=[S'(\lambda):D'(\mu)]$
whenever $D'(\mu)\ne0$.

In latter sections we will be particularly interested in the analogues
of the Young modules in this setup. For each $\lambda\in\Lambda^+$ let
$P'(\lambda)$ be the projective cover of $L'(\lambda)$. Suppose that
$\mu\in\Lambda^+$ and let $\mathcal N(\mu)=\phi_\mu'\Schur'(\Lambda)$,
where $\phi_\mu'$ is the identity map on $N(\mu)$. Then $\phi_\mu'$
is an idempotent so $\mathcal N(\mu)$ is a projective
$\Schur'(\Lambda)$--module. Therefore, there exist non--negative 
integers $c_{\lambda\mu}\ge0$ such that 
$$\mathcal N(\mu)\cong P'(\mu)
         \oplus\bigoplus_{\lambda\gdom\mu}c_{\lambda\mu}P'(\lambda).$$ 
In fact, because $\mathcal N(\mu)_\ZZ=\mathcal M(\mu)_\ZZ'$ it follows
by a specialization argument that the integers $c_{\lambda\mu}$ are
the same as those appearing in Proposition~\ref{M^mu S-filtration}.  We call
$Y'(\mu)=\Sfun(P'(\lambda))$ a {\sf twisted Young module}.

Define the filtration multiplicities $[Y'(\lambda):S'(\mu)]$ exactly
as in section~3. 

\begin{Proposition}
\label{twisted Youngs}
Suppose that $R$ is a field and let $\lambda$ and $\mu$ be 
multipartitions of $n$. Then
\begin{enumerate}

\item $N(\mu)\cong
Y'(\mu)\oplus\bigoplus_{\lambda\gdom\mu}Y'(\nu)^{c_\lambda\mu}$
where the integers $c_{\lambda\mu}$ are the same as those appearing 
in Proposition~\ref{M^mu S-filtration};
\item $Y'(\mu)$ is indecomposable;
\item $Y'(\lambda)\cong Y'(\mu)$ if and only if $\lambda=\mu$;
\item the Young module $Y'(\mu)$ has a dual Specht filtration in which
the number of subquotients equal to $S'(\lambda)$ is
$[\W'(\lambda):L'(\mu)]$; and,
\item $[Y'(\mu):S'(\lambda)]=[\W'(\lambda):L'(\mu)].$
\end{enumerate}
\end{Proposition}

\begin{proof}
This can be proved in exactly the same way as in
Theorem~\ref{Young modules}; alternatively, one can use a specialization
argument. 
\end{proof}

We remark that the set of Young modules
$\set{Y(\lambda)|\lambda\vdash n}$ and the set of twisted Young
modules $\set {Y'(\lambda)|\lambda\vdash n}$ do not usually coincide;
however, we always have that
$$\set{Y(\lambda)|D(\lambda)\ne0}=\set{Y'(\lambda)|D'(\lambda)\ne0},$$
because these modules are the indecomposable direct summands of $\H$.
T see this use Corollary~\ref{M^mu summands} to show that if $\lambda$ is
a multipartition and $D(\lambda)\ne0$ then $Y(\lambda)$ is the
projective cover of $D(\lambda)$ and, similarly, that
$Y'(\mu)$ is the projective cover of $D'(\mu)$. It follows that
$D(\lambda)\cong D'(\mu)$ if and only if $Y(\lambda)\cong Y'(\mu)$. By
the results of \cite{FLOTW,Ariki:can,DM:Morita,A:class}, the
correspondence between these two different labellings of the simple
$\H$--modules is given by a generalization of Kleshchev's version of
the Mullineux map (that is, in terms of paths in the associated
crystal graphs).

\section{Contragredient duality}

We now investigate contragredient duality for the category of
$\H$--modules; this will give us the connection between the
Specht modules with the dual Specht modules constructed in the
previous section. The aim of the section is really to construct a dual
Specht filtration of $M(\lambda)$; in essence, this is the
main tool that we need in order to understand the tilting modules
of the cyclotomic Schur algebras.

Recall that $*$ is the unique anti--isomorphism of $\H$ such that
$T_i^*=T_i$ for $0\le i<n$. Given a right $\H$--module $M$ define its
{\sf contragredient dual} $M^\dual$ to be the dual module
$\Hom_R(M,R)$ equipped with the right $\H$--action $(\phi
h)(m)=\phi(mh^*)$ for all $\phi\in M^\dual$, $h\in \H$ and $m\in M$.
A module $M$ is {\sf self--dual} if $M\cong M^\dual$. By standard
arguments, $M$ is self--dual if and only if $M$ possesses a
non--degenerate associative bilinear form (the form $\<\ ,\ \>$ is
associative if $\<xh,y\>=\<x,yh^*\>$ for all $x,y\in M$ and $h\in\H$).

If $M$ is a submodule of $\H$ the reader should be careful not to
confuse the dual module $M^\dual$ with
$M^*=\set{m^*|m\in M}$.

Constructing dual bases inside $\H$ is, in general, quite hard. We are
going to do it by comparing the two bases $\{m_{\s\t}\}$ and
$\{n_{\s\t}\}$ of $\H$. First we need to introduce some notation for
conjugate multipartitions and tableaux.

Recall that the conjugate of a composition $\sigma$ is the partition
$\sigma'=(\sigma_1',\sigma_2',\dots)$ where $\sigma_i'$ is the number
of nodes in column $i$ of the diagram of $\sigma$.  If
$\lambda=\rtuple\lambda$ is a multicomposition then the {\sf conjugate} of
$\lambda$ is the multipartition 
$\lamp=\((\lambda^{(r)})',\dots,(\lambda^{(1)})'\)$. Observe that if
$\lambda\gedom\mu$ then $\mu'\gedom\lambda'$ (and conversely, if
$\lambda$ and $\mu$ are multipartitions).

If $\T=\rtuple\T$ is a $\lambda$--tableau of type $\mu$ then the
{\sf conjugate} of~$\T$ is the $\lamp$--tableau $\T'=\rtuple{\T'}$ where
$\T'(i,j,s)=\T(j,i,r-s+1)$ for all $(i,j,s)\in[\lambda']$; that is,
${\T'}^{(s)}$ is the tableau obtained by interchanging the rows and
columns of $\T^{(r-s+1)}$. Notice that $\Type(\T')=\Type(\T)$.
Further, if $\t$ is a tableau of type $\omega$ then $\t$ is standard
if and only if $\t'$ is standard. 

For each multipartition $\lambda$ let $\tllam=(\t^\lamp)'$; thus,
$\tllam$ is the standard $\lambda$--tableau with the numbers
$1,2,\dots,n$ entered in order first down the columns of
$\t_\lambda^{(r)}$ and then the columns of $\t_\lambda^{(r-1)}$ and so
on. Observe that if $\t$ is a standard $\lambda$--tableau then
$\tlam\gedom\t\gedom\tllam$. We set $w_\lambda=d(\tllam)$. The
following Lemma is well--known; it can be proved by induction on~$\t$.

\begin{Lemma}
Suppose that $\lambda$ is a multipartition of $n$ and that $\t$
is a standard $\lambda$--tableau. Then $d(\t)d(\t')^{-1}=w_\lambda$
and $\len(w_\lambda)=\len(d(\t))+\len(d(\t'))$.  
\label{top}\end{Lemma}

We also extend the dominance order to pairs of tableaux in the usual
way.

If $\t$ is a standard tableau and $k$ an integer with $1\le k\le n$
then the {\sf residue} of $k$ in $\t$ is defined to be
$\res_\t(k)=q^{2(j-i)}Q_s$ if $k$ appears in row $i$ and column $j$ of
component $s$ of $\t$. Residues are important because of the following
result.

\begin{Point}
{\it}{(James--Mathas\cite[Prop.~3.7]{JM:cyc-Schaper})}
Suppose that $\s$ and $\t$ are standard $\lambda$--tableau and that
$1\le k\le n$. Then there exist $r_{\u\v}\in R$ such that
$$m_{\s\t}L_k=\res_\t(k)m_{\s\t}
         +\sum_{(\u,\v)\gdom(\s,\t)}r_{\u\v}m_{\u\v}.$$
\label{L_k action}\end{Point}

Let $\K=\Q(q,Q_1,\dots,Q_r)$ and, following \cite{M:gendeg}, define
$F_\t\in\HK$ by
$$F_\t=\prod_{k=1}^n
       \prod_{\substack{c\in\mathcal R(k)\\c\ne\res_\t(k)}}
       \frac{L_k-c}{\res_\t(k)-c},$$
where $\mathcal R(k)
   =\set{q^{2d}Q_s|1\le s\le r\And|d|<k\And d\ne0\If r=1\And k=2,3}$.
Finally, given two standard $\lambda$--tableau $\s$ and $\t$ set
$f_{\s\t}=F_\s m_{\s\t}F_\t$ and $g_{\s\t}=F_{\s'}n_{\s\t}F_{\t'}$.
From the definitions, $\(\res_\t(k)\)'=\res_{\t'}(k)$ in $\ZZ$, for
all tableau $\t$ and all $k$. This implies that $F_\t'=F_{\t'}$ and
hence that $g_{\s\t}=f_{\s\t}'$ in $\HK$; see \cite{M:gendeg}.
 
Using (\ref{L_k action}) we obtain the following.

\begin{Point}
{\it}{(Mathas~\cite{M:gendeg})}\label{ortho props} Suppose that $\H=\HK$.
\begin{enumerate}
\item 
$m_{\s\t}=f_{\s\t}+\sum_{\a,\b}r_{\a\b}f_{\a\b}$
for some $r_{\a\b}\in R$ with $r_{\a\b}\ne0$ only if
$(\a,\b)\gdom(\s,\t)$.
\item 
$n_{\s\t}=g_{\s\t}+\sum_{\a,\b}r_{\a\b}g_{\a\b}$
for some $r_{\a\b}\in R$ with $r_{\a\b}\ne0$ only if
$(\a,\b)\gdom(\s,\t)$.
\item Suppose that $\s,\t,\u$ and $\v$ are standard tableaux. Then
$f_{\s\t}g_{\u\v}=0$ unless $\t=\u'$. 
\end{enumerate}\end{Point}   

By (i) and (ii), the sets $\{f_{\s\t}\}$ and $\{g_{\s\t}\}$ are both
bases of $\HK$. In fact, by \cite{M:gendeg}, both bases are
self--orthogonal with respect to the bilinear form $\<\ ,\ \>$ which
we introduce below.

\begin{Lemma}
Suppose that $\s$ and $\t$ are standard $\lambda$--tableaux and
that $\v$ and $\u$ are standard $\mu$--tableaux such that
$m_{\s\t}n_{\v\u}\ne0$. Then $\v'\gedom\t$.
\label{killer}
\end{Lemma}

\begin{proof}
Now, $m_{\s\t}n_{\v\u}\ne0$ in $\H$ only if
$m_{\s\t}n_{\v\u}\ne0$ in $\HZ$ since $\H\cong\HZ\otimes_\Z R$ and
specialization maps the standard basis of $\HZ$ to the standard basis
of $\H$ and, similarly, for the $n_{\v\u}$ basis elements. Hence, by
embedding $\HZ$ into $\HK$ in the natural way, we
may assume that $\H=\HK$.  By parts (i) and (ii) of (\ref{ortho props})
there exists scalars $r_{\a\b},r_{\c\d}\in R$ such that
$$0\ne m_{\s\t}n_{\v\u}
   =\Big(f_{\s\t}+\sum_{(\a,\b)\gdom(\s,\t)}r_{\a\b}f_{\a\b}\Big)
    \Big(g_{\v\u}+\sum_{(\c,\d)\gdom(\v,\u)}r_{\c\b}f_{\c\d}\Big).$$
Therefore, there exist $(\a,\b)\gedom(\s,\t)$ and
$(\c,\d)\gedom(\v,\u)$ such that $f_{\a\b}g_{\c\d}\ne0$; so,
$\c'=\b$ by (\ref{ortho props})(iii). Consequently,
$\v'\gedom\c'=\b\gedom\t$ as required.
\end{proof}

It is possible to give a direct proof of Lemma~\ref{killer} without using
the two orthogonal bases $\{f_{\s\t}\}$ and $\{g_{\s\t}\}$ of $\H$
(cf.~\cite[Lemma~4.11]{murphy:basis}); however, the proof above is
both easier and nicer because it avoids long calculations with
the relations in~$\H$. (The proof of (\ref{ortho props}) is
straightforward and also avoids such calculations.)

Recall from (\ref{AK basis}) that
$\set{L_1^{a_1}\dots L_n^{a_n}T_w|0\le a_i<r\And w\in\Sym_n}$
is a basis of $\H$. Define
$\tau\map\H R$ to be the $R$--linear map determined by
$$\tau(L_1^{a_1}\dots L_n^{a_n} T_w)=\begin{cases}
           1,&\If a_1=\dots=a_n=0\And w=1,\\
           0,&\Otherwise.
\end{cases}$$
This map was introduced by Bremke and Malle~\cite{BM:redwds}
who showed that $\tau$ is a trace form; that is, $\tau(ab)=\tau(ba)$
for all $a,b\in\H$. (The definition above is slightly different from
Bremke and Malle's; it is shown in \cite{MM:trace} that the two
definitions coincide.) Combining the definition with the fact that
$\tau$ is a trace form shows that $\tau(h^*)=\tau(h)$, for all
$h\in\H$.

Define a bilinear form $\<\ ,\ \>\map{\H\times\H}\H$ on $\H$ by
$\<h_1,h_2\>=\tau(h_1h_2^*)$ for all $h_1,h_2\in\H$. Then $\<\ ,\ \>$
is an associative bilinear form on $\H$; further, $\<\ ,\ \>$ is
symmetric because $\tau$ is a trace form.

For each multipartition $\lambda$ set
$$Q_\lambda=(-1)^{n(r-1)}\Prod_{s=1}^rQ_s^{n-|\lambda^{(s)}|}.$$
Then $Q_\lambda\in R$ and $Q_\lambda$ is a unit if and
only if $Q_s$ is a unit whenever $|\lambda^{(s)}|<n$.

Many of the results which follow rely upon the following result.

\begin{Theorem}
Suppose that $(\s,\t)$ is a pair of $\lambda$--tableau and that
$(\u,\v)$ are $\mu$--tableaux. Then
$$\<m_{\s\t},n_{\u\v}\>=
    \begin{cases} Q_\lambda, &\If (\u',\v')=(\s,\t),\\
                   0, &\If (\u',\v')\not\gedom(\s,\t).
\end{cases}$$
\label{orthogonal}
\end{Theorem}

\begin{proof}
Suppose first that $\<m_{\s\t},n_{\u\v}\>\ne0$. Now
$\<m_{\s\t},n_{\u\v}\>=\tau(m_{\s\t}n_{\v\u})$,
so $m_{\s\t}n_{\v\u}\ne0$; hence, $\v'\gedom\t$ by Lemma~\ref{killer}.
Now $\tau$ is a trace form and $\tau(h)=\tau(h^*)$,
for all $h\in\H$; so, applying these two facts, we have
$\tau(m_{\s\t}n_{\v\u})=\tau(n_{\v\u}m_{\s\t})=\tau(m_{\t\s}n_{\u\v})$; 
hence, $m_{\t\s}n_{\u\v}\ne0$ and $\u'\gedom\s$ by Lemma~\ref{killer}.
Therefore, if $\<m_{\s\t},n_{\u\v}\>\ne0$ then $(\u',\v')\gedom(\s,\t)$.

Now assume that $(\u',\v')=(\s,\t)$. Then 
$T_{w_\lambda}=T_{d(\t)}T_{d(\t')}^*=T_{d(\s')}T_{d(\s)}^*$ by 
Lemma~\ref{top}. Therefore, once again using the fact that $\tau$ is a trace
form,
\begin{align*}
\<m_{\s\t},n_{\s'\t'}\>
   &=\tau\(m_{\s\t}n_{\t'\s'}\)
    =\tau\(T_{d(\s)}^*m_\lambda T_{d(\t)}
          T_{d(\t')}^*n_\lamp T_{d(\s')}\)\\
   &=\tau\(T_{d(\s')}T_{d(\s)}^*m_\lambda T_{w_\lambda}n_\lamp)
    =\tau\(T_{w_\lambda}^*m_\lambda T_{w_\lambda} n_\lamp).
\end{align*}    
Finally, 
$\tau(T_{w_\lambda}^*m_\lambda T_{w_\lambda} n_\lamp\)=Q_\lambda$ 
by \cite[Prop.~5.12]{M:gendeg}, so we're done.
\end{proof}

As a first consequence we obtain a new proof that $\H$ is a symmetric
algebra.

\begin{Corollary}
[Malle--Mathas \cite{MM:trace}] 
Suppose that $q,Q_1,\dots,Q_r$ are invertible elements of $R$. Then
$\<\ ,\ \>$ is a non--degenerate associative symmetric bilinear form
on $\H$. Therefore, $\H$ is a symmetric algebra; in particular, it is
self--dual
\label{trace}\end{Corollary} 

At first sight, this proof of Corollary~\ref{trace} is considerably easier than
the original proof in~\cite{MM:trace}; however, all of the work is
hidden in the calculation of $\tau(T_{w_\lambda}^*m_\lambda
T_{w_\lambda} n_\lamp\)$ from \cite{M:gendeg} and this is quite
involved. The payoff for this extra effort is Theorem~\ref{orthogonal} which
shows that the two bases $\{m_{\s\t}\}$ and $\{n_{\s\t}\}$ are almost
orthogonal; this fact will be used many times in what follows.

Let $\lambda$ be a multipartition. We next show that 
$S'(\lambda)\cong S(\lambda')^\dual$, and so justify the term 
{\it dual} Specht module.  Recall that $N(\lambda)=n_\lambda\H$ and
that $S'(\lambda)$ is a quotient of $N(\lambda)$. 

\begin{Corollary}
Suppose that $q,Q_1,\dots,Q_r$ are invertible elements of $R$ and
let $\lambda$ be a multipartition of $n$.  Then 
$S'(\lambda)\cong S(\lambda')^\dual$.
\label{dual Spechts}
\end{Corollary}

\begin{proof}
Now, $S(\lambda')$ is a submodule of $\H/\H(\lamp)$
and $S'(\lambda)$ is a submodule of~$\H/\Hlambarp$. By
Theorem~\ref{orthogonal}(i) the modules $\H(\lamp)$ and $\Hlambarp$ are
orthogonal with respect to the form $\<\ ,\ \>$, as are $M(\lamp)$ and
$\Hlambarp$, and $N(\lambda)$ and $\H(\lamp)$.  Therefore, $\<\ ,\ \>$
induces an associative bilinear form 
$\<\ ,\ \>_{S(\lambda)}\map{S(\lamp)\times S'(\lambda)}R$ given by
$$\< a+\Hlambar,b+\Hlambarp\>_{S(\lambda)}=\<a,b\>=\tau(ab^*).$$
In particular, if $\s\in\Std(\lamp)$ and $\t\in\Std(\lambda)$ then
$$\< m_\s,n_\t\>_{S(\lambda)}
       =\begin{cases} Q_\lambda, &\If \t'=\s,\\
                        0,&\Unless \t'\gedom\s,
\end{cases}$$
by Theorem~\ref{orthogonal}.  Hence, $\<\ ,\ \>_{S(\lambda)}$ is non--degenerate
and $S'(\lambda)\cong S(\lamp)^\dual$ as required.  
\end{proof}

Recall from (\ref{sstd basis thm}) that $M(\lambda)$ is free as an
$R$--module with basis
$$\set{m_{\S\t}|\S\in\SStd(\lambda,\mu)\And\t\in\Std(\lambda)
                   \For\lambda\vdash n}.$$
It was shown in \cite[Cor.~4.15]{DJM:cyc} that this basis gives rise
to a Specht filtration of $M(\mu)$. Similarly, the basis of 
Corollary~\ref{N sstd basis thm} produces a dual Specht filtration of
$N(\lambda)$.  We next produce another basis of $M(\mu)$ which
exhibits a dual Specht filtration of $M(\mu)$ and, similarly, a basis
of $N(\mu)$ which exhibits a Specht filtration of $N(\mu)$. As a
byproduct we will also obtain a non--degenerate associative bilinear
form on each of these modules and hence see that they are both
self--dual. 

A $\lambda$--tableau $\T$ is {\sf column semistandard} if
$\T'$ is semistandard. If $\mu\in\Lambda$ and $\lambda\in\Lambda^+$ let
$$\CStd(\lambda,\mu)=\set{\T|\T'\in\SStd(\lambda,\mu)}$$ 
be the set of column semistandard $\lambda$--tableaux of type $\mu$.
Observe that if $\CStd(\lambda,\mu)\ne\0$ then 
$\SStd(\lambda',\mu)\ne\0$, so $\lambda'\gedom\mu$; equivalently, 
$\mu'\gedom\lambda$. We also set
$\CStd(\Lambda^+,\mu)
      =\bigcup_{\lambda\in\Lambda^+}\CStd(\lambda,\mu)$.

As a final piece of notation, if $\v$ is any tableau and $1\le k\le n$
then write $\comp_\v(k)=s$ if $k$ appears in component~$s$ of $\v$. 

\begin{Lemma}
Suppose that $\mu$ is a multicomposition and that 
$m_\mu n_{\u\v}\ne0$ or $n_\mu m_{\u\v}\ne0$ for some standard
tableaux $\u$ and $\v$. Then $\mu(\u)$ is column semistandard. 
\label{semikiller}
\end{Lemma}

\begin{proof}
As in the proof of Lemma~\ref{killer} we may assume that $\H=\HZ$. We
consider only the case where $m_\mu n_{\u\v}\ne0$; the other case can
be proved by applying the involution $'$. 

Before we begin the proof proper we remark that it is well--known, and
easy enough to check, that if $s_i\in\Sym_\mu$ then $m_\mu T_i=qm_\mu$;
similarly, if~$s_j\in\Sym_\lambda$ then 
$T_j n_\lambda=-q^{-1}n_\lambda$. 

First, because $\u'$ is standard the entries of $\mu(\u')$ are weakly
increasing along rows. Suppose that the entries of $\mu(\u')$ are not
strictly increasing down columns. Then we can find integers $i<j$ such
that $i$ and $j$ are in the same row of $\t^\mu$ and the same column
of $\u'$. The entries in $\t^\mu$ are consecutive so this means that
there is exists an integer $i$ such that $i$ and $i+1$ are in the same
row of $\t^\mu$ and the same row of $\u$. Therefore,
$$qm_\mu n_{\u\v}=(m_\mu T_i)n_{\u\v}
                      =m_\mu(T_i n_{\u\v})
                      =-q^{-1}m_\mu n_{\u\v}.$$
Consequently, $m_\mu n_{\u\v}=0$ since $\HZ$ is $\ZZ$-free. 

It remains to show that $\mu(\u')$ satisfies condition (iii) of
Definition~\ref{semistandard}. If $\mu$ is a multipartition then
$m_\mu=m_{\t^\mu\t^\mu}$; so $m_{\t^\mu\t^\mu}n_{\u\v}\ne0$ and
$\u'\gedom\t^\mu$ by Lemma~\ref{killer}. Looking at the definitions,
we see that $\mu(\u')$ satisfies condition Definition~\ref{semistandard}(iii)
because $\u'\gedom\t^\mu$. Hence,~$\u$ is column semistandard as
claimed.

If $\mu$ is a multicomposition (and not a multipartition) let
$\vec\mu=\rtuple{\vec\mu}$ be the multipartition obtained by ordering
the parts in each component $\mu^{(s)}$ of $\mu$. Then we can find a
permutation~$d$ of minimal length in
$\Sym_{|\mu|}=\Sym_{|\mu^{(1)}|}\times\dots\times\Sym_{|\mu^{(r)}|}$
such that $d\Sym_\mu=\Sym_{\vec\mu}d$. Then $T_dm_\mu=m_{\vec\mu}
T_d$.  Now, $m_\mu n_{\u\v}$ is non--zero so $T_dm_\mu
n_{\u\v}=m_{\vec\mu}T_dn_{\u\v}$ is also non--zero. Therefore, there
exists tableaux $\a$ and $\b$ such that $n_{\a\b}$ is a non--zero
summand of $T_d n_{\u\v}$ and $m_{\vec\mu}n_{\a\b}\ne0$.  By the last
paragraph $\vec\mu(\a)$ satisfies Definition~\ref{semistandard}(iii).  This
implies that $\mu(\u)$ also satisfies Definition~\ref{semistandard}(iii) because
$\comp_\u(k)=\comp_\a(k)$, for $1\le k\le n$, by
\cite[Prop.~3.18]{DJM:cyc} since $d\in\Sym_{|\mu|}$. Hence, $\mu(\u)$
is column semistandard.  
\end{proof}

If $\S$ is a $\lambda$--tableau of type $\mu$ let $\dot\S$ be the
unique standard tableau such that $\mu(\dot\S)=\S$ 
and~$\dot\S\gedom\s$ whenever $\s$ is a standard $\lambda$--tableau with
$\mu(\s)=\S$. The tableau $\dot\S$ is denoted
$\operatorname{first}(\S)$ in \cite{JM:cyc-Schaper}. The permutation
$d(\dot\S)$ is a distinguished $(\Sym_\lambda,\Sym_\mu)$--double coset
representative; that is, it is the unique element of minimal length in
$\Sym_\lambda d(\dot\S)\Sym_\mu$. We emphasize
that~$\dot\S$ is a {\it standard} tableau.

\begin{Proposition}
Suppose that $\mu$ is a multicomposition of $n$.
Then $M(\mu)$ is free as an $R$--module with basis
$\set{m_\mu n_{\dot\S\t}|\S\in\CStd(\lambda,\mu)\And\t\in\Std(\lambda)
                           \ForSome\lambda\vdash n}$
and $N(\mu)$ is free as an $R$--module with basis
$\set{n_\mu m_{\dot\S\t}|\S\in\CStd(\lambda,\mu)\And\t\in\Std(\lambda)
                           \ForSome\lambda\vdash n}$.
\label{dual basis}
\end{Proposition}

\begin{proof}
We only prove for the claim for $M(\mu)$; the second statement
can be proved by a similar argument, or by specialization.

By (\ref{ortho props})(iii) $\{n_{\s\t}\}$ is a basis of $\H$, so
$M(\mu)$ is spanned by the elements $m_\mu n_{\s\t}$, where $\s$ and
$\t$ range over all pairs of standard tableaux of the same shape.
Furthermore, if $m_\mu n_{\s\t}\ne0$ then $\mu(\s)$ is column
semistandard by Lemma~\ref{semikiller}. Hence, $M(\mu)$ is spanned by the
elements $m_\mu n_{\s\t}$ with $\mu(\s)$ column semistandard. Now, if
$d(\s)$ and $d(\u)$ are in the same $(\Sym_\lambda,\Sym_\mu)$--double
coset then $m_\mu T_{d(\s)}n_\lambda=\pm q^a m_\mu T_{d(\u)}n_\lambda$
for some integer~$a$; see the remarks at the start of the proof of
Lemma~\ref{semikiller}. By definition $d(\dot\S)$ is the unique element of
minimal length in its double coset; therefore, the elements in the
statement of the Lemma span $M(\mu)$. However, now we are done because
$M(\mu)$ is $R$--free and the number of elements in our spanning set
is exactly the rank of $M(\mu)$ by (\ref{sstd basis thm}).
\end{proof}

Combining Lemma~\ref{semikiller} and the Proposition we have.

\begin{Corollary}
Suppose that $\mu$ is a multicomposition and that $\s$ and $\t$
are standard tableaux. Then $m_\mu n_{\s\t}\ne0$ if and only if
$\mu(\s)$ is column semistandard. Similarly, $n_\mu m_{\s\t}\ne0$ if
and only if $\mu(\s)$ is column semistandard. 
\label{killed}\end{Corollary}

Using Proposition~\ref{dual basis}, the argument of Lemma~\ref{Weyl filtration} produces
the following result.

\begin{Corollary}
\label{dual Specht series}
Suppose that $\mu$ is a multicomposition of $n$. Then there exist
filtrations
$$M(\mu)=M_1\supset\dots\supset M_k\supset M_{k+1}=0\quad\And\quad
  N(\mu)=N_1\supset\dots\supset N_k\supset N_{k+1}=0$$ 
of $M(\mu)$ and $N(\mu)$, respectively, and multipartitions
$\lambda_1,\dots,\lambda_k$, such that $\mu'\gedom\lambda_i$,
$M_i/M_{i+1}\cong S'(\lambda_i)$
and $N_i/N_{i+1}\cong S(\lambda_i),$ for $1\le i\le k$.  Moreover, if
$\lambda$ is any multipartition of $n$ then $\#\set{1\le i\le
k|\lambda_i=\lambda}=\#\CStd(\lambda,\mu)$.
\end{Corollary}

\begin{Remark}
As an $R$--module, $N_i$ is the submodule of $N(\mu)$ with
basis the set of elements $n_\mu m_{\dot\S\t}$ with
$\Shape(\S)\gedom\lambda_i$, for $1\le i\le k$. In particular,
$S(\mu')\cong N_k$ is spanned by $\set{n_\mu
m_{\t_{\mu'}\t|\t\in\Std(\mu')}}$; note that $\mu(\t_{\mu'})$ is the
unique column standard $\mu'$--tableau of type~$\mu$. Therefore,
$S(\mu')\cong n_\mu T_{w_{\mu'}}^* m_{\mu'}\H$; this is a result of Du
and Rui~\cite{DuRui:branching}. Similarly, $S'(\mu')\cong m_\mu
T_{w_{\mu'}}^* n_{\mu'}\H$.
\label{remark}\end{Remark}

Because $S(\lambda')\cong S(\lambda')^\dual$, the Specht filtrations
of $M(\mu)$ given by \cite[Cor.~4.15]{DJM:cyc} and the last result
suggest that $M(\lambda)$ is self--dual. A similar remark applies to
$N(\lambda)$. When $r=1$ it is clear that both of these modules are
self--dual because they are induced representations from parabolic
subalgebras. 

We need a non--degenerate associative bilinear form.
Let $\<\ ,\ \>_\mu$ be the bilinear map on~$M(\mu)$ determined by
$$\<m_{\S\t}, m_\mu n_{\dot\U\v}\>_\mu
        =\<m_{\S\t},n_{\dot\U\v}\>,$$
where $m_{\S\t}$ and $m_\mu n_{\dot\U\v}$ run over  the bases of
(\ref{sstd basis thm}) and Proposition~\ref{dual basis}, respectively.

If $\S$ is a semistandard tableau let $\S_{(i,s)}$ be the subtableau
of $\S$ consisting of those entries $(j,t)$ with $(j,t)\preceq(i,s)$
(see Definition~\ref{semistandard}). We extend the dominance order to the
set of semistandard tableaux by defining $\S\gedom\T$ if
$\Shape(\S_{(i,s)})\gedom\Shape(\T_{(i,s)})$ for all $(i,s)$. This
definition coincides with our previous definition of dominance when
$\S$ is a standard tableau (recall that we are identifying standard
tableaux and semistandard tableaux of type~$\omega$).

\begin{Proposition}
Suppose that  $Q_1,\dots,Q_r$ are invertible elements of $R$ and
that $\mu$ is a multicomposition. Then $\<\ ,\ \>_\mu$ is a
non--degenerate associative bilinear form on $M(\lambda)$. In
particular, $M(\mu)$ is self--dual. Similarly, $N(\mu)$ is self--dual.
\label{MN self-dual}
\end{Proposition}

\begin{proof}
We prove the Proposition only for $M(\mu)$; the result for
$N(\mu)$ can be obtained using a similar argument or by
specialization.

Suppose that $\S\in\SStd(\lambda,\mu)$, $\t\in\Std(\lambda)$,
$\U\in\CStd(\rho,\mu)$ and $\v\in\Std(\rho)$ for some 
multipartitions~$\lambda$ and $\rho$. Applying the definitions we 
find that
$$\<m_{\S\t},m_\mu n_{\dot\U\v}\>_\mu
     =\<m_{\S\t},n_{\dot\U\v}\>
      =\sum_{\substack{\s\in\Std(\lambda)\\\mu(\s)=\S}}
           \<m_{\s\t},n_{\dot\U\v}\>.$$
Therefore, by Theorem~\ref{orthogonal}, $\<m_{\S\t},m_\mu n_{\dot\U\v}\>_\mu=0$
unless there is a standard tableau $\s$ such that
$(\dot\U',\v')\gedom(\s,\t)$ and $\mu(\s)=\S$; hence,
$(\U',\v')\gedom(\S,\t)$. (Here
$\dot\U'=(\dot\U)'$ and not~$(\U')^{\mbox{\normalsize$\cdot$}}$; in 
general these tableaux are different.)  

Next suppose that $(\U',\v')=(\S,\t)$.
Then $\dot\U=(\S')^\cdot\gedom\s'$ whenever $\s'$ is a standard tableau
with $\mu(\s')=\S'$; therefore, $\s\gedom\dot\U'$ whenever $\s$ is a
standard tableau such that $\mu(\s)=\S$. Therefore, if
$(\U,\v)=(\S',\t')$ then
$$\<m_{\S\t},n_{\dot\U\v}\>_\mu
     =\sum_{\substack{\s\in\Std(\lambda)\\\mu(\s)=\S}}
           \<m_{\s\t},n_{\dot\U\v}\>
     =\<m_{\dot\U'\v'},n_{\dot\U\v}\>
     =Q_\lambda$$ 
by Theorem~\ref{orthogonal}. 

Combining the last two paragraphs shows that the matrix
$\(\<m_{\S\t},m_\mu n_{\dot\U\v}\>_\mu\)$ is invertible and, hence,
that the form $\<\ ,\ \>_\mu$ on $M(\mu)$ is non--degenerate.

The harder part is to prove that $\<\ ,\ \>_\mu$ is associative. Now,
the form $\<\ ,\ \>$ on $\H$ is associative, so if~$h\in\H$ then
$$\<m_{\S\t}h,m_\mu n_{\dot\U\v}\>_\mu
    =\<m_{\S\t}h, n_{\dot U\v}\>
    =\<m_{\S\t},n_{\u\v}h^*\>
    =\sum_{\a,\b}r_{\a\b}\<m_{\S\t},n_{\a\b}\>,
$$
where $n_{\a\b}h^*=\sum_{\a,\b}r_{\a\b}n_{\a\b}$ for some 
$r_{\a\b}\in R$. Write $m_{\S\t}=m_\mu h_{\S\t}^\mu$, for some 
$h_{\S\t}^\mu\in\H$. Then
$$\<m_{\S\t},n_{\a\b}\>=\tau(m_{\S\t}n_{\b\a})
                       =\tau(n_{\b\a}m_{\S\t})
                       =\tau(n_{\b\a}m_\mu h_{\S\t}^\mu)
                       =\tau(h_{\t\S}^\mu m_\mu n_{\a\b}),
$$
where the second equality uses the fact that $\tau$ is a trace form
and the last equality follows because $\tau(a)=\tau(a^*)$ for all
$a\in\H$.  Therefore, if $\<m_{\S\t},n_{\a\b}\>\ne0$ then $\mu(\a)$ is
column semistandard by Lemma~\ref{semikiller}. Let $\A=\mu(\a)$ and recall
that $m_\mu n_{\a\b}=\pm q^i m_\mu n_{\dot\A\b}$ for some integer $i$
(which depends on $\a$); write $m_\mu n_{\a\b}=\lambda_\a m_\mu
n_{\dot\A\b}$ and set $r_{\A\b}=\sum_\a \lambda_ar_{\a\b}$, where the sum
runs over those standard tableaux with $\mu(\a)=\A$. Then we have
shown that
$$\<m_{\S\t}h,m_\mu n_{\dot\U\v}\>_\mu
     =\sum_{\substack{\b\in\Std(\lambda)\\
                       \A\in\CStd(\lambda,\mu)}}
                       r_{\A\b}\<m_{\S\t},n_{\dot\A\b}\>.
$$
On the other hand, by Lemma~\ref{semikiller} again,
\begin{align*}
\<m_{\S\t},m_\mu n_{\dot\U\v}h^*\>_\mu
     &=\sum_{\a,\b}r_{\a\b}\<m_{\S\t}, m_\mu n_{\a\b}\>_\mu
     =\sum_{\substack{\a,\b\in\Std(\lambda)\\
                       \mu(\a)\in\CStd(\lambda,\mu)}}
         r_{\a\b}\<m_{\S\t},m_\mu n_{\a\b}\>_\mu\\
     &=\sum_{\substack{\b\in\Std(\lambda)\\
                       \A\in\CStd(\lambda,\mu)}}
         r_{\A\b}\<m_{\S\t},m_\mu n_{\dot\A\b}\>_\mu
      =\sum_{\substack{\b\in\Std(\lambda)\\
                       \A\in\CStd(\lambda,\mu)}}
         r_{\A\b}\<m_{\S\t},n_{\dot\A\b}\>.
\end{align*}
Hence, 
$\<m_{\S\t}h,m_\mu n_{\dot\U\v}\>_\mu
  =\<m_{\S\t},m_\mu n_{\dot\U\v}h^*\>_\mu$, so the form is
associative as claimed.
\end{proof}

\begin{Corollary}
Suppose that $R$ is a field and let $\lambda$ be a multipartition
of $n$. Then both the Young module $Y(\lambda)$ and the twisted Young
module $Y'(\lambda)$ are self--dual.
\label{self dual Youngs}
\end{Corollary}

\begin{proof}
By the Proposition, $M(\lambda)$ and $N(\lambda)$ are both
self--dual. Hence, the result follows by induction on $\lambda$ using
Theorem~\ref{Young modules} and Proposition~\ref{twisted Youngs}, respectively.
\end{proof}

\section{The cyclotomic tilting modules}

Let $(A,X^+)$ be a quasi--hereditary algebra, where $X^+$ is the poset of
weights for $A$; see, for example, \cite[Appendix]{Donkin:book}. For
each $\lambda\in X^+$ there is a standard module $\Delta(\lambda)$ with
simple head $\L(\lambda)$ and a costandard module with simple socle
$\L(\lambda)$. An $A$--module $M$ has a $\Delta$--filtration if it has
a filtration in which every subquotient isomorphic to a standard
module; similarly, $M$ has a $\nabla$--filtration if every subquotient
is isomorphic to a costandard module. An $A$--module $T$ is a {\sf tilting
module} if it has both a $\Delta$--filtration and a
$\nabla$--filtration.

\begin{Point}
{\it}{(Ringel~\cite{Ringel})} Suppose that $R$ is a field and that
$(A,X^+)$ is a quasi--hereditary algebra. Then, for each 
$\lambda\in X^+$, there is a unique indecomposable tilting module
$T(\lambda)$ such that
$$[T(\lambda):\Delta(\lambda)]=1\quad\And\quad
  [T(\lambda):\Delta(\mu)]\ne0\OnlyIf\lambda\ge\mu.$$ 
Moreover, if $T$ is any tilting module then 
$$T\cong \bigoplus_{\lambda\in X^+} T(\lambda)^{t_\lambda}$$ 
for some non--negative integers~$t_\lambda$.
\label{Ringel}\end{Point}

The $T(\lambda)$ are the {\sf partial tilting modules} of $A$. A {\sf full
tilting module} for $A$ is any tilting module which contains every
$T(\lambda)$, for $\lambda\in X^+$, as a direct summand.

By \cite[Cor.~6.18]{DJM:cyc} the cyclotomic Schur algebras are
quasi--hereditary algebras with weight poset $\Lambda^+$. The standard
modules of $\Schur(\Lambda)$ are the Weyl modules and the costandard
modules are their contragredient duals. In this section we will
describe the partial tilting modules of $\Schur(\Lambda)$ when
$\omega\in\Lambda$ and the parameters $Q_1,\dots,Q_r$ are distinct and
non--zero.

First consider the case $r=1$. Suppose that $d\ge1$ and let
$\Lambda_{d,n}$ be the set of compositions of $n$ into at most $d$
parts and let $V$ be a free $R$--module of rank $d$. Then $\H(\Sym_n)$
acts on $V^{\otimes n}$ (by $q$--analogues of place permutations) and
$V^{\otimes n}\cong M(\Lambda)$; see \cite{DJ:qWeyl}. The
Dipper--James~\cite{DJ:Schur} $q$--Schur algebra $\Schur(d,n)$ is the
cyclotomic $q$--Schur algebra $\Schur(\Lambda_{d,n})$; by the above
remarks $\Schur(d,n)\cong\End_{\H(\Sym_n)}(V^{\otimes n})$.
Donkin~\cite{Donkin:tilt} has shown that when $d\ge n$ the tilting
modules for $\Schur(d,n)$ are the indecomposable direct summands of
the exterior powers $\wedge^\lambda V
=\wedge^{\lambda_1}V\otimes\dots\otimes\wedge^{\lambda_d}V$.  In the
cyclotomic case we do not have a description
of~$M(\Lambda)=\bigoplus_\mu M(\mu)$ as a tensor product;
nevertheless, we do have the following analogue of the exterior
powers.

\begin{Definition}
For each multicomposition $\alpha$ let
$E(\alpha)
    =\Hom_\H\(M(\Lambda),N(\alpha)\)$.
\end{Definition}

By definition, $E(\alpha)$ is a right $\Schur(\Lambda)$--module.
Recall that $\vec\alpha$ is the multipartition obtained by reordering
the parts of $\alpha^{(s)}$ for each~$s$. By the argument of Corollary~\ref{M^mu summands}, $N(\alpha)\cong N(\vec\alpha)$; therefore, $E(\alpha)\cong
E(\vec\alpha)$.  Hence, there is no loss in assuming that $\alpha$ is
a multipartition.  The $E(\alpha)$ are very similar to the modules
$\M(\mu)=\Hom_\H\(M(\Lambda),M(\mu)\)$ of Proposition~\ref{M^mu S-basis}.
The $E(\alpha)$ play the role of exterior powers and the 
$\M(\mu)$ the symmetric powers.

We will show that $E(\alpha)$ has a Weyl filtration and that it is
self--dual; hence, it also has a dual Weyl filtration. This will
enable us to show that the tilting modules of $\Schur(\Lambda)$ are
the indecomposable summands of the $E(\lambda)$ as $\lambda$ runs over
the multipartitions in $\Lambda^+$.

The next result is a first step towards producing a basis for
$E(\alpha)$. By general principles, if $f\in N(\alpha)\cap M(\mu)^*$
then left multiplication by $f$ is an $\H$--module homomorphism
from $M(\mu)$ into $N(\alpha)$; in fact, every element of $E(\alpha)$
arises in this way.

\begin{Lemma}
Suppose that $\alpha$ and $\mu$ are multicompositions of $n$. 
Then there is an isomorphism of $R$--modules
$\Hom_\H\(M(\mu),N(\alpha)\)\cong N(\alpha)\cap M(\mu)^*$
given by $\theta\mapsto\theta(m_\mu)$. In particular, if
$\theta\in\Hom_\H\(M(\mu),N(\alpha)\)$ and $h_\theta=\theta(m_\mu)$
then $h_\theta\in N(\alpha)\cap M(\mu)^*$ and 
$\theta(m)=h_\theta m$, for all $m\in M(\mu)$.
\label{annihilator}
\end{Lemma}

\begin{proof}
This follows from Theorem~5.16 and Lemma~5.2 of \cite{DJM:cyc}.
(Theorem~5.16 says that the double annihilator,
$\set{h\in\H|hs=0\Whenever m_\mu s=0}$, of $m_\mu$ is $\H m_\mu$;
Lemma~5.2 observes that this property of the double annihilator implies 
the Lemma.)
\end{proof}

Therefore, to give a basis of $E(\alpha)$ it is enough to find a basis
of $N(\alpha)\cap M(\mu)^*$. To do this we need to make the following
assumption.

\begin{Assumption}
For the rest of this paper assume that 
$Q_1,\dots,Q_r$ are distinct.
\label{assumption}\end{Assumption}

The only place where we explicitly use Assumption~\ref{assumption} is in the
proof of the following theorem; however, almost everything which
follows relies on this result. Unless otherwise stated, this
assumption will remain in force for the rest of the paper.

\begin{Theorem}
Suppose that $Q_1,\dots,Q_r$ are all distinct and let $\alpha$
and $\mu$ be multicompositions of $n$.  Then $N(\alpha)\cap M(\mu)^*$
is free as an $R$--module with basis
$$\set{n_\alpha m_{\dot\S\T}|\S\in\CStd(\lambda,\alpha),
              \T\in\SStd(\lambda,\mu)\ForSome\lambda\vdash n}.$$
\label{NM basis}
\end{Theorem}

\begin{proof}
First note that $N(\alpha)\cap M(\mu)^*$ is free because it is
a submodule of a free module. Next, by (\ref{sstd basis thm}), if
$\S\in\CStd(\lambda,\alpha)$ and $\T\in\SStd(\lambda,\mu)$, for some
multipartition $\lambda$, then 
$n_\alpha m_{\dot\S\T}\in N(\alpha)\cap M(\mu)^*$ since $N(\alpha)$ is
a right ideal and $M(\mu)^*$ is a left ideal of $\H$. Moreover, as
$\S$ and $\T$ run over the possible choices above, these elements
remain linearly independent because the elements $n_\alpha
m_{\dot\S\t}$ are a basis of $N(\alpha)$ by Proposition~\ref{dual basis}. 

It remains to show that these elements span $N(\alpha)\cap M(\mu)^*$.
It will be convenient to let~$\ASStd(\Lambda^+,\mu)$ be the set of
tableaux of type $\mu$ which satisfy conditions (i) and (ii) of
Definition~\ref{semistandard} but {\it not} condition (iii) --- these are
``almost'' semistandard tableaux. 

By Proposition~\ref{dual basis}, if $x\in N(\alpha)\cap M(\mu)^*$ then 
$x=\sum r_{\U\v}n_\alpha m_{\dot\U\v}$, for some $r_{\U\v}\in R$,
where the sum is over all pairs $(\U,\v)$ of tableaux with
$\U$ column semistandard of type~$\alpha$ and $\v$ standard.  Now, if
$(i,i+1)\in\Sym_\mu$ then $m_\mu T_i=q$; hence, $xT_i=qx$ and, as
in \cite[(4.19)]{murphy:basis} (compare \cite[Lemma~4.11]{DJM:cyc}),
it follows that if $r_{\U\v}\ne0$ then $i$ and $i+1$ are not in the
same column of~$\v$ and that $r_{\S\t}=r_{\S\v}$ where $\t=\v(i,i+1)$;
that is, $\V=\mu(\v)$ satisfies conditions (i) and (ii) of
Definition~\ref{semistandard}.  Therefore, $x=\sum r_{\U\V} n_\alpha
m_{\dot\U\V}$ where the sum is over pairs $(\U,\V)$ with
$\U\in\CStd(\Lambda^+,\alpha)$ and
$\V\in\SStd(\Lambda^+,\mu)\cup\ASStd(\Lambda^+,\mu)$. By the first
paragraph, $n_\alpha m_{\dot\S\T}\in N(\alpha)\cap M(\mu)^*$ when $\T$
is semistandard, so we may assume that $r_{\U\V}=0$ unless
$\V\in\ASStd(\Lambda^+,\mu)$. Thus, we are reduced to showing that if
we have an element $x\in N(\alpha)\cap M(\mu)^*$ which can be written
in the form 
$$x=\sum_{\substack{\U\in\CStd(\Lambda^+,\alpha)\\
                    \V\in\ASStd(\lambda,\mu)}}
           r_{\U\V}n_\alpha m_{\dot\U\V},$$
for some $r_{\U\V}\in R$, then $x=0$. By way of contradiction, suppose
that $x\ne0$.

Fix $(\S,\T)$ with $r_{\S\T}\ne0$ such that $r_{\U\V}=0$ whenever
$(\S,\T)\gdom(\U,\V)$ for $\U\in\CStd(\Lambda^+,\alpha)$ and
$\V\in\ASStd(\Lambda^+,\mu)$. Let $\ddot\T$ be the unique standard
tableau such that $\mu(\ddot\T)=\T$ and $\t\gedom\ddot\T$ whenever
$\mu(\t)=\T$; the tableau $\ddot\T$ is denoted
$\operatorname{last}(T)$ in \cite{JM:cyc-Schaper}. Let $i$ be the
smallest positive integer such that
$c=\comp_{\ddot\T}(i)>\comp_{\t^\mu}(i)$; such an $i$ exists because
$\T$ does not satisfy Definition~\ref{semistandard}(iii). If $j<i$ then
$\comp_{\ddot\T}(j)<\comp_{\ddot\T}(i)$ by the minimality of $i$ and
the fact that $\comp_{\t^\mu}(j)\le\comp_{\t^\mu}(i)$; in particular,
this implies that~$i$ must appear in the first row and first column of
${\ddot\T}^{(c)}$.  Following \cite{DJM:cyc} define
$$y_i=T_{i-1}\dots T_1\prod_{s=1}^{\comp_{\t^\mu}(i)}(L_1-Q_s).$$ 
Then $y_i\ne0$ since $\comp_{\t^\mu}(i)<c\le r$. Moreover, $m_\mu y_i=0$ 
by \cite[Lemma~5.8]{DJM:cyc} (when translating into the notation of
\cite{DJM:cyc} note that $\gamma_i=\comp_{\t^\mu}(i)$). Let
$\t=\ddot\T s_{i-1}\dots s_1$; then
$\res_{\t}(1)=\res_{\ddot\T}(i)=Q_c$ and
$\len(d(\t))=\len(d(\ddot\T))+i-1$; so, 
$m_{\dot\S\ddot\T}T_{i-1}\dots T_1=m_{\dot\S\t}$. Furthermore, by the
cancellation property of the Bruhat--Chevalley order,
$\v\gedom\ddot\T$ if and only if $\v s_{i-1}\dots s_1\gedom\t$ since
$\len(d(\t))=\len(d(\ddot\T))+i-1$; see, for example,
\cite[Cor.~3.9]{M:ULect}. Therefore, using (\ref{L_k action}) for the
third equality, we have
\begin{align*}
x y_i&=\Big(r_{\S\T}n_\alpha m_{\dot\S\ddot\T} 
       +\sum_{\substack{\u,\v\\\ddot\T\not\gedom\v}}
             r_{\u\v} n_\alpha m_{\u\v}\Big) y_i\\
   &=\Big(r_{\S\T}n_\alpha m_{\dot\S\t}
       +\sum_{\substack{\u,\v\\\t\not\gedom\v}}
             r_{\u\v}'n_\alpha m_{\u\v}\Big) 
         \prod_{s=1}^{\comp_{\t^\mu}(i)}(L_1-Q_s)\\
   &=r_{\S\T}\prod_{s=1}^{\comp_{\t^\mu}(i)}(Q_c-Q_s)\cdot
      n_\alpha m_{\dot\S\t}
 +\sum_{\substack{\u,\v\\\t\not\gedom\v}}
             r_{\u\v}''n_\alpha m_{\u\v} y_i
\end{align*}
for some $r_{\u\v}',r_{\u\v}''\in R$. By Assumption~\ref{assumption} the
coefficient of $n_\alpha m_{\dot\S\t}$ in $x y_i$ is non--zero because
$c>\comp_{\t^\mu}(i)$ --- and $r_{\S\T}\ne0$; note also that 
$n_\alpha m_{\dot\S\t}\ne0$ by Corollary~\ref{killed}. Therefore, $xy_i\ne0$.
However, $x\in N(\alpha)\cap M(\mu)^*$ and, as remarked above, 
$m_\mu y_i=0$, so this contradicts the assumption that $x\ne0$. 
Consequently, $x=0$ and the Theorem follows.
\end{proof}

\begin{Remark}
If $Q_s=Q_t$, where $s\ne t$, then the rank of
$n_\alpha\H\cap\H m_\mu$ can be larger than that predicted by 
Theorem~\ref{NM basis}. For example, suppose that $Q_1=Q_2$, when $r=n=2$, and 
take $\alpha=\mu=((1),(1))$. Then $n_\alpha=m_\mu=L_1-Q_2$ and, by
(\ref{S-basis}),
$$n_\alpha\H\cap\H m_\mu=M(\mu)\cap M(\mu)^*$$
is the free $R$--module with basis
$\set{m_{\S\T}|\S,\T\in\SStd(\Lambda^+,\mu)}$; this is
an $R$--module of rank~$3$. In contrast, if $Q_1\ne Q_2$ then 
by Theorem~\ref{NM basis}
$$n_\alpha\H\cap\H m_\mu=Rn_\alpha T_1m_\mu
                        =R(L_1-Q_1)T_1(L_1-Q_2);$$
this time the intersection has rank $1$. (By direct a calculation,
$(L_1-Q_1)T_1(L_1-Q_2)$ is an element of $M(\mu)\cap M(\mu)^*$ if and 
only if $Q_1=Q_2$.)

Shoji has shown that if the parameters $Q_1,\dots,Q_r$ are distinct
then $M(\mu)$ is an induced module (more accurately, he has shown that
there exists an induced module which has the same image as $M(\mu)$ in
the Grothendieck group of $\H$). It should be possible to prove
analogue of Frobenius reciprocity for the modules $M(\mu)$ and
$N(\alpha)$ using Shoji's work; this would give a better explanation
as to why the rank of $\Hom_\H\(M(\mu),N(\alpha)\)$ is independent of
the choice of parameters $Q_1,\dots,Q_r$ in the presence of
Assumption~\ref{assumption}.
\end{Remark}

Now we reap some consequences of Theorem~\ref{NM basis}. We emphasize that
even though we do not explicitly state Assumption~\ref{assumption} it remains in
force for all of these results.

\begin{Corollary}
Suppose that $\alpha$ and $\mu$ are multicompositions of $n$.
Then 
$$N(\alpha)\cap M(\mu)^*=n_\alpha\H m_\mu.$$
\label{intersection}
\end{Corollary}

\begin{proof}
Certainly, $n_\alpha\H m_\mu\ss N(\alpha)\cap M(\mu)^*$.
Conversely, by Theorem~\ref{NM basis} a basis of $N(\alpha)\cap M(\mu)^*$ is
given by the elements $n_\alpha m_{\dot\S\T}$, where
$\S\in\CStd(\Lambda^+,\alpha)$ and $\T\in\SStd(\Lambda^+,\mu)$.
As $n_\alpha m_{\dot\S\T}\in n_\alpha\H m_\mu$, by 
(\ref{sstd basis thm}), we also have the opposite inclusion.
\end{proof}

Notice that Theorem~\ref{NM basis} and Corollary~\ref{intersection} imply that if 
$\mu$ is a multipartition then
$$n_\mu\H\cap\H m_{\mu'}
     =n_\mu\H m_{\mu'}
         =Rn_\mu m_{\t_{\mu'}\t^{\mu'}}
     =Rn_\mu T_{w_\mu} m_{\mu'}.$$
Du and Rui~\cite{DuRui:branching} have shown that 
$n_\mu\H m_{\mu'}=Rn_\mu T_{w_\mu} m_{\mu'}$
(without assuming that the parameters~$Q_s$ are distinct).
This is interesting because the element 
$n_\mu T_{w_\mu} m_{\mu'}$ generates the Specht module
$S(\mu')$; see Remark~\ref{remark}.

\begin{Corollary}
Suppose that $\Hom_\H\(M(\mu),N(\alpha)\)\ne0$ for some
multicompositions $\alpha$ and $\mu$. Then $\alpha'\gedom\mu$.
\end{Corollary}

\begin{proof}
By Lemma~\ref{annihilator} $N(\alpha)\cap M(\mu)^*\ne0$ since
$\Hom_\H\(M(\mu),N(\alpha)\)\ne0$. Therefore, we can find a
multipartition $\lambda$ and tableaux $\S\in\CStd(\lambda,\alpha)$ and
$\T\in\SStd(\lambda,\mu)$ such that 
$n_\alpha m_{\dot\S\T}$ is a non--zero element of
$N(\alpha)\cap M(\mu)^*$. Hence,
$\alpha'\gedom\lambda\gedom\mu$, so $\alpha'\gedom\mu$ as required.
\end{proof}

Similarly, if $\Hom_\H\(N(\alpha),M(\mu)\)\ne0$ then
$\mu'\gedom\alpha$.

\begin{Corollary}
Suppose that $\alpha$ and $\mu$ are multicompositions of $n$.
Then
$$\set{n_{\S\dot\T}m_\mu|\S\in\SStd(\lambda,\alpha),
              \T\in\CStd(\lambda,\mu)\ForSome\lambda\vdash n}$$
is a basis of $N(\alpha)\cap M(\mu)^*$.
\label{MN basis}
\end{Corollary}

\begin{proof}
By Theorem~\ref{NM basis} the $R$--module $N(\alpha)\cap M(\mu)^*$ is
stable under specialization (or, if you prefer, base change);
therefore, it is enough to consider the case $\H=\HZ$.  Applying the
involutions $'$ and $*$ to the basis of $N(\mu)\cap M(\alpha)^*$ given
by Theorem~\ref{NM basis} yields the result. Alternatively, this can be proved
by modifying the argument of Theorem~\ref{NM basis}.
\end{proof}

These results allow us to give two bases for $E(\alpha)$.

\begin{Definition}
Suppose that $\lambda$ is a multipartition of $n$ and that
$\alpha$ and $\mu$ are two multicompositions of $n$. For 
tableaux
$\S\in\CStd(\lambda,\alpha)$, $\T\in\SStd(\lambda,\mu)$,
$\A\in\SStd(\lambda,\alpha)$ and $\B\in\CStd(\lambda,\mu)$
let $\theta_{\S\T}$ and~$\theta_{\S\T}'$
be the homomorphisms in $E(\alpha)$ determined by
$$\theta_{\S\T}(m_\nu h)=\delta_{\nu\mu}n_\alpha m_{\dot\S\T}h
\quad\And\quad
\theta_{\A\B}'(m_\nu h)=\delta_{\nu\mu}n_{\A\dot\B}m_\mu h,$$
for all $h\in\H$ and all $\nu\in\Lambda$.
\end{Definition}

Lemma~\ref{annihilator} together with Theorem~\ref{NM basis} and Corollary~\ref{MN basis},
respectively, show that these maps are elements of~$E(\alpha)$.
Indeed, these results show that each of the corresponding sets of such
maps is a basis of $E(\alpha)$.  More precisely, we have the
following.

\begin{Proposition}
Suppose that $Q_1,\dots,Q_r$ are all distinct and let $\alpha$
be a multipartition of~$n$. Then $E(\alpha)$ is free as an $R$--module
with bases $\mathscr E$ and $\mathscr E'$ where
\begin{align*}
\mathscr E&=\set{\theta_{\S\T}|\S\in\CStd(\lambda,\alpha)
                \And\T\in\SStd(\lambda,\mu)\ForSome\lambda\vdash n}\\
\text{and\ }
\mathscr E'&=\set{\theta_{\A\B}'|\A\in\SStd(\lambda,\alpha)
    \And\B\in\CStd(\lambda,\mu)\ForSome\lambda\vdash n}.
\end{align*}
\label{E basis}\end{Proposition}

Now that we have the required notation it is a good time to note that
$E(\alpha)$ is cyclic.

\begin{Corollary}
Suppose that $\omega\in\Lambda$ and that $\alpha$ is a
multicomposition. Then $E(\alpha)$ is a cyclic
$\Schur(\Lambda)$--module; more precisely,
$E(\alpha)=\theta_{\T^\omega_\alpha\T^\omega}\Schur(\Lambda)$ where
$\T^\omega_\alpha=\alpha(\t^\omega)$.
\label{E cyclic}
\end{Corollary}

\begin{proof}
The map $\theta_{\T^\omega_\alpha\T^\omega}$ is the extension
to $E(\alpha)$ of the homomorphism $\H\To N(\alpha)$ given by
$\theta_{\T^\omega_\alpha\T^\omega}(h)
              =n_\alpha m_{\dot\T^\omega_\alpha\T^\omega}h
              =n_\alpha h$,
for all $h\in\H$.  Suppose that in $\S\in\CStd(\lambda,\alpha)$ and
$\T\in\SStd(\lambda,\alpha)$ for some multipartition $\lambda$.Then
$\theta_{\S\T}=\theta_{\T^\omega_\alpha\T^\omega}\phi_{\dot\S\T}$
(both maps send $m_\mu$ to~$n_\alpha m_{\dot\S\T}$), so $E(\alpha)$ is
cyclic as claimed.
\end{proof}

As an application of Proposition~\ref{E basis} we now show that
$E(\alpha)$ has a Weyl filtration.

\begin{Theorem}
Suppose that $\omega\in\Lambda$ and that $Q_1,\dots,Q_r$ are
all distinct. Then there exist
multipartitions $\lambda_1,\dots,\lambda_k$ in $\Lambda$ and an
$\Schur(\Lambda)$--module filtration
$$E(\alpha)= E_1\supset\dots\supset E_k\supset E_{k+1}=0$$
such that$E_i/ E_{i+1}\cong\W(\lambda_i)$ and
$\alpha'\gedom\lambda_i$, for $i=1,\dots,k$. Moreover, if $\lambda$ is
any multipartition in $\Lambda^+$ then
$\#\set{1\le i\le k|\lambda_i=\lambda}=\#\CStd(\lambda,\alpha)$.
\label{E filtration}
\end{Theorem}

\begin{proof}
For the most part this is a familiar argument; however, towards
the end there is a small twist so we give the details. Let
$\S_1,\dots,\S_k$ be the complete set of column semistandard tableaux
of type $\alpha$ ordered so that $i>j$ whenever
$\Shape(\S_i)\gedom\Shape(\S_j)$. For each $i$ let
let~$\lambda_i=\Shape(\S_i)$; then $\alpha'\gedom\lambda_i$ since
$\CStd(\lambda_i,\alpha)\ne\0$.

For $i=1,\dots,k$ let $E_i$ be the $R$--submodule of $E(\alpha)$
spanned by the elements $\theta_{\S_j\T}$ with $j\ge i$ and
$\T\in\SStd(\lambda_j,\Lambda)$. To prove the Theorem it is enough to
show that, $E_i$ is a submodule of $E(\alpha)$ and that
$E_i/E_{i+1}\cong\W(\lambda_i)$ for each $i$.

Suppose that $i\ge1$ and let $\T\in\SStd(\lambda_i,\mu)$,
$\U\in\SStd(\rho,\sigma)$ and $\V\in\SStd(\rho,\nu)$, for some
multicompositions $\mu,\nu,\sigma\in\Lambda$ and $\rho\in\Lambda^+$.
Consider the product $\theta_{\S_i\T}\phi_{\U\V}$. By definition,
$\theta_{\S_i\T}\phi_{\U\V}=0$ unless
$\mu=\Type(\T)=\Type(\U)=\sigma$; so suppose that $\sigma=\mu$. In
order to write $\theta_{\S_i\T}\phi_{\U\V}$ as a linear combination of
the basis elements of $E(\alpha)$ it suffices to consider
$\(\theta_{\S_i\T}\phi_{\U\V}\)(m_\nu)$. Write 
$m_{\U\V}=m_\mu h^\mu_{\U\V}$, for some $h^\mu_{\U\V}\in\H$. Then we 
have
\begin{align*}
\(\theta_{\S_i\T}\phi_{\U\V}\)(m_\nu)
     &=\theta_{\S_i\T}(m_{\U\V})
      =\theta_{\S_i\T}(m_\mu)h^\mu_{\U\V}
      =n_\alpha m_{\dot\S_i\T}h^\mu_{\U\V}\\
     &\equiv n_\alpha\Big(\sum_{\Y\in\SStd(\lambda_i,\nu)}
                    r_{\S_i\Y}m_{\dot\S_i\Y}\Big)
              \mod{\mathcal H^{\lambda_i}}\\
     &=\sum_{\Y\in\SStd(\lambda_i,\nu)}r_{\S_i\Y}n_\alpha m_{\dot\S_i\Y},
\intertext{%
  where the second line follows from (\ref{W mult}). Hence,
  Theorem~\ref{NM basis} implies that
}
\theta_{\S_i\T}\phi_{\U\V}&
   \equiv\sum_{\Y\in\SStd(\lambda_i,\nu)}r_{\S_i\Y}\theta_{\S_i\Y}
             \mod{E_{i+1}}.
\end{align*}
All of our claims now follow.
\end{proof}

\begin{Remark}
By Theorem~\ref{E filtration} the module $\W(\lambda)$ is a Weyl
module composition factor of~$E(\alpha)$ whenever
$\CStd(\lambda,\alpha)\ne\0$; thus, $\alpha'\gedom\lambda$. Moreover,
since $\CStd(\omega,\alpha)$ is always non--empty this means that
$\W(\omega)$ is always a composition factor of $E(\alpha)$;
consequently, the assumption that $\omega\in\Lambda$ is necessary in
Theorem~\ref{E filtration}.  If $\omega\notin\Lambda$ then $E(\alpha)$ is
still an $\Schur(\Lambda)$--module (for all multicompositions
$\alpha$); however, we are not able to give a Weyl filtration of
$E(\alpha)$ in this case.
\end{Remark}

\begin{Corollary}
Suppose that $R$ is a field, $\omega\in\Lambda$ and that 
$\lambda$ and $\mu$ are multipartitions of~$n$. Then 
$[E(\lambda){:}\W(\lambda')]=1$ and
$[E(\lambda){:}\W(\mu)]\ne0$ only if $\lambda'\gedom\mu$.
\label{multiplicities}
\end{Corollary}

\begin{proof}
By the Theorem, $[E(\lambda){:}\W(\mu)]=\#\CStd(\mu,\lambda)$.
Therefore, if $[(\lambda){:}\W(\mu)]\ne0$ then
$\CStd(\mu,\lambda)\ne\0$ and $\mu'\gedom\lambda$. Finally,
$[E(\lambda){:}\W(\lambda')]=1$ because
$\lambda(\t_{\lambda'})=(\Tlam)'$ is the unique column semistandard
$\lambda'$--tableau of type $\lambda$ (just as
$\T^\lambda=\lambda(\t^\lambda)$ is the unique semistandard
$\lambda$--tableau of type $\lambda$). 
\end{proof}

As in the previous section, if $E$ is an $\Schur(\Lambda)$--module
then its {\sf contragredient dual} $E^\dual$ is the dual space
$\Hom_R(E,R)$ equipped with the contragredient action:
$(f\phi)(x)=f(x\phi^*)$, for $f\in E^\dual$,
$\phi\in\Schur(\Lambda)$ and $x\in E$. Again, $E$ is {\sf self--dual} if
$E\cong E^\dual$.

If $E$ is an $\Schur(\Lambda)$--module and $\mu\in\Lambda$ let
$E_\mu=E\phi_\mu$ be the {\sf $\mu$--weight space} of $E$.

\begin{Theorem}
Suppose that $Q_1,\dots,Q_r$ are distinct invertible elements
of $R$ and let $\alpha$ be a multicomposition. Then $E(\alpha)$ is
self--dual.
\label{E dual}
\end{Theorem}

\begin{proof}
Define a bilinear map 
$\{\ ,\ \}_\alpha\map{E(\alpha)\times E(\alpha)}R$ by
$$\{\theta_{\S\T},\theta_{\A\B}'\}_\alpha
      =\begin{cases} \<m_{\dot\S\T},n_{\A\dot\B}\>,&\If\Type(\T)=\Type(\B),\\
               0,&\Otherwise.
\end{cases}$$
where $\theta_{\S\T}$ and $\theta_{\A\B}'$ run over the two bases
$\mathscr E$ and $\mathscr E'$ of $E(\alpha)$ from Proposition~\ref{E basis}. 
By definition, the different weight spaces
$E_\mu$, $\mu\in\Lambda$, of $E(\alpha)$ are
orthogonal with respect to  $\{\ ,\ \}_\alpha$. Suppose
then that $\Type(\T)=\Type(\B)$ and let $\lambda=\Shape(\S)$. Then, as
in the proof of Proposition~\ref{MN self-dual}, 
$$\{\theta_{\S\T},\theta_{\A\B}'\}_\alpha
      =\begin{cases} Q_\lambda,&\If (\A',\B')=(\S,\T),\\
              0,       &\If(\A',\B')\not\gedom(\S,\T),
\end{cases}$$
where $\lambda=\Shape(\S)$. Hence, $\{\ ,\ \}_\alpha$ is a
non--degenerate bilinear form on $E(\alpha)$.

Once again, the harder task is to prove that $\{\ ,\ \}_\alpha$ is
associative. Choose tableaux $\S,\T,\A$ and $\B$ as above. Suppose that
$\phi\in\Hom_\H(M(\nu),M(\mu))$ for some $\nu,\mu\in\Lambda$. Then
$\theta_{\S\T}\phi\in E_\nu$ and 
$\theta_{\A\B}'\phi^*\in E_\mu$; therefore, if $\mu\ne\nu$ then 
$$\{\theta_{\S\T}\phi,\theta_{\A\B}\}_\alpha=0
     =\{\theta_{\S\T},\theta_{\A\B}\phi^*\}_\alpha$$
because different weight spaces are orthogonal with respect to $\{\ ,\
\}_\alpha$. Suppose then that $\mu=\nu$. By considering weight spaces
we may also assume that $\Type(\T)=\mu=\Type(\B)$. Write
$\phi(m_\mu)=m_\mu h$ for some $h\in\H$; then
$(\theta_{\S\T}\phi)(m_\mu h)=n_\alpha m_{\dot\S\T}h$, so
$\theta_{\S\T}\phi$ is determined by $n_\alpha m_{\dot\S\T}h$.
Similarly, since $\phi(m_\mu)\in m_\mu\H\cap\H m_\mu$ we can also
write $\phi(m_\mu)=\tilde hm_\mu$. Therefore,
$\phi^*(m_\mu)=m_\mu{\tilde h}^*=(\tilde hm_\mu)^*=(m_\mu h)^*=h^*m_\mu$;
consequently,
$(\theta_{\A\B}'\phi*)(m_\mu)=\theta_{\A\B}'(m_\mu)\tilde h^*
                             =n_{\A\dot\B}m_\mu{\tilde h}^*
                             =n_{\A\dot\B}h^*m_\mu.$
The bilinear form $\<\ ,\ \>$ on
$\H$ is associative so, as in the proof of Proposition~\ref{MN self-dual}, 
we have
$$\{\theta_{\S\T}\phi,\theta_{\A\B}\}_\alpha
        =\<m_{\dot\S\T}h,n_{\A\dot\B}\>
        =\<m_{\dot\S\T},n_{\A\dot\B}h^*\>
        =\{\theta_{\S\T},\theta_{\A\B}\phi^*\}_\alpha.$$
Therefore, $\{\ ,\ \}_\alpha$ is associative and the proof is complete.
\end{proof}

Combining the last two results we obtain our main theorem.

\begin{Theorem}
\label{tilting modules}
Suppose that $R$ is a field, $\omega\in\Lambda$ and
that $Q_1,\dots,Q_r$ are distinct non--zero elements of $R$.
\begin{enumerate}
\item If $\lambda\in\Lambda^+$ then
$$E(\lambda)\cong 
  T(\lambda')\oplus\bigoplus_{\lambda'\gdom\mu}T(\mu)^{e_{\lambda\mu}}$$
for some non--negative integers $e_{\lambda\mu}$.
\item The tilting modules of
$\Schur(\Lambda)$ are the indecomposable direct summands of the
modules $\set{E(\lambda)|\lambda\in\Lambda^+}.$
\end{enumerate}
\end{Theorem}

\begin{proof}
By Theorem~\ref{E filtration}, $E(\lambda)$ has a $\Delta$--filtration;
therefore, $E(\lambda)$ also has a $\nabla$--filtration since
$E(\lambda)$ is self--dual by Theorem~\ref{E dual}. Hence, $E(\lambda)$ is a
tilting module. Furthermore, by 
Corollary~\ref{multiplicities}, $[E(\lambda){:}\W(\lambda')]=1$ and
if $[E(\lambda){:}\W(\mu)]>0$ then $\lambda'\gedom\mu$. Therefore,
by Ringel's theorem~(\ref{Ringel}), there exist non--negative 
integers $e_{\lambda\mu}$ such that
$$E(\lambda)\cong 
  T(\lambda')\oplus\bigoplus_{\lambda'\gdom\mu}T(\mu)^{e_{\lambda\mu}}.$$
This proves (i). Part (ii) now follows by induction on the dominance
order using (\ref{Ringel}).
\end{proof}

\section{Ringel duality}

We now turn our attention to the Ringel dual of $\Schur(\Lambda)$. By
definition, the {\sf Ringel dual} of~$\Schur(\Lambda)$ is the algebra
$\End_{\Schur(\Lambda)}(T)$, where $T$ is any full tilting module for
$\Schur(\Lambda)$; thus, the Ringel dual is determined only up to
Morita equivalence. By Theorem~\ref{tilting modules} the module
$$E(\Lambda)=\bigoplus_{\alpha\in\Lambda}E(\alpha)
       =\bigoplus_{\alpha,\mu\in\Lambda}\Hom_\H\(M(\mu),N(\alpha)\)$$
is a full tilting module for $\Schur(\Lambda)$ when $\omega\in\Lambda$
and the parameters $Q_1,\dots,Q_r$ are distinct and non--zero.

Let $\alpha$ and $\beta$ be two multicompositions. Then for any
$\S\in\SStd(\lambda,\beta)$ and $\T\in\SStd(\lambda,\alpha)$ there is
an $\H$--module homomorphism $\phiST'\map{N(\alpha)}N(\beta)$; this
induces an $\Schur(\Lambda)$--module homomorphism
$\Phi_{\S\T}'\map{E(\alpha)}E(\beta)$ given by
$\Phi_{\S\T}'(\theta)=\phiST'\theta$, for $\theta\in E(\alpha)$.  In
fact, we will show that these give all of the
$\Schur(\Lambda)$--module homomorphisms from $E(\alpha)$ to
$E(\beta)$.

The next two results do not require that the parameters 
$Q_1,\dots,Q_r$ be distinct.

\begin{Proposition}
Suppose that $\omega\in\Lambda$ and let $\alpha$ and $\beta$
be multicompositions of $n$. Then
$\Hom_{\Schur(\Lambda)}\(E(\alpha),E(\beta)\)$ is free as an
$R$--module with basis
$$\set{\Phi_{\S\T}'|\S\in\SStd(\lambda,\beta),
                   \T\in\SStd(\lambda,\alpha)
                   \ForSome\lambda\in\Lambda^+}.$$
\label{endo}
\end{Proposition}

\begin{proof}
As indicated above the maps $\Phi_{\S\T}'$ belong to 
$\Hom_{\Schur(\Lambda)}\(E(\alpha),E(\beta)\)$. Moreover, they are
linearly independent because the $\phi_{\S\T}'$ are a basis
of $\Hom_\H\(N(\alpha),N(\beta)\)$ by Proposition~\ref{S' basis}(i). Thus, it
remains to see that these maps span 
$\Hom_{\Schur(\Lambda)}\(E(\alpha),E(\beta)\)$.

Suppose that $\Phi\in\Hom_{\Schur(\Lambda)}\(E(\alpha),E(\beta)\)$.
If $\theta\in\Hom_\H\(M(\mu),N(\alpha)\)$ then $\Phi(\theta)$ belongs
to $\Hom_\H\(M(\mu),N(\beta)\)\cong N(\alpha)$; that is, $\Phi$ maps weight spaces to
weight spaces. Now, $E(\alpha)$ is generated by
$\theta_{\T^\omega_\alpha\T^\omega}$ by Corollary~\ref{E cyclic}; so, $\Phi$ is
determined by $\Phi(\theta_{\T^\omega_\alpha\T^\omega})$.  Moreover,
$\Phi(\theta_{\T^\omega_\alpha\T^\omega})
         \in\Hom_\H\(\H,N(\beta)\)\cong N(\beta)$
since $\Phi$ maps weight spaces to weight spaces. Therefore, 
$\Phi(\theta_{\T^\omega_\alpha\T^\omega})
         \in\Hom_\H\(N(\alpha),N(\beta)\)$
since $\Phi$ is an $\Schur(\Lambda)$--module homomorphism and
$\H\cong\Hom_\H(\H,\H)$ (where we identify $h\in\H$ with left
multiplication by $h$). Hence, we can write 
$\Phi(\theta_{\T^\omega_\alpha\T^\omega})
         =\sum_{\S,\T}r_{\S\T}\phi_{\S\T}'$
for some $r_{\S\T}\in R$ by Proposition~\ref{S' basis}(i). Therefore, 
$\Phi=\sum_{\S,\T}r_{\S\T}\Phi_{\S\T}'$, completing the proof.
\end{proof}

If $\Lambda$ is a saturated set of multicompositions let
$$E(\Lambda)=\bigoplus_{\alpha\in\Lambda}E(\alpha)
      =\Hom_\H\(M(\Lambda),N(\Lambda)\)
      =\bigoplus_{\mu,\alpha\in\Lambda}\Hom_\H\(M(\mu),N(\alpha)\)
.$$
Then $E(\Lambda)$ is an
$\(\Schur'(\Lambda),\Schur(\Lambda)\)$--bimodule. Moreover, it has the
following double centralizer property.

If $A$ is an algebra let $A^{\text{op}}$ be the opposite algebra in
which the order of multiplication is reversed.

\begin{Corollary}
Suppose that $\omega\in\Lambda$.
Then there are canonical isomorphisms of $R$--algebras
$$\End_{\Schur(\Lambda)}\(E(\Lambda)\)
           \cong\Schur'(\Lambda)^{\text{op}}
\quad\And\quad
  \End_{\Schur(\Lambda)'}\(E(\Lambda)\)\cong\Schur(\Lambda).$$
\end{Corollary}

\begin{proof}
The first isomorphism, $\Phi_{\S\T}'\mapsto\phi_{\S\T}'$,
is immediate from (the proof of) Proposition~\ref{endo}. The second isomorphism
follows by symmetry.
\end{proof}

As a special case we have a description of the Ringel dual of
$\Schur(\Lambda)$.

\begin{Corollary}
Suppose that $R$ is a field, $\omega\in\Lambda$ and that
$Q_1,\dots,Q_r$ are distinct invertible elements of $R$.
Then the Ringel dual of $\Schur(\Lambda)$ is isomorphic to
$\Schur'(\Lambda)^{\text{op}}$.
\end{Corollary}

Finally, we want to determine the $\nabla$--filtration multiplicities
in the tilting modules. We actually don't need to do any work here
because the general theory of tilting modules for quasi--hereditary
algebras tells us that
$[T(\lambda):\nabla(\mu)]=[\W(\mu'),L(\lambda')]$ (see, for example,
\cite[Appendix]{Donkin:book}); however, we want to show how this
result can be derived using Young modules and Specht filtrations. 

Recall that $\Sfun\map{\Schur(\Lambda)\Mod}\H\Mod; M\mapsto M\phio$ is
the Schur functor.

\begin{Proposition}
Suppose that $R$ is a field, $\omega\in\Lambda$ and let
$\lambda$ be a multipartition. Then
$\Sfun\(T(\lambda')\)\cong Y'(\lambda)$ as $\H$--modules.
\label{tilting Youngs}
\end{Proposition}

\begin{proof}
Applying the definitions
$\Sfun\(E(\lambda)\)\cong\Hom_\H\(\H,N(\lambda)\)\cong N(\lambda)$,
where the last isomorphism comes from Lemma~\ref{annihilator} (or directly).
Hence, the Schur functor $\Sfun$ induces an injective map from
$\End_{\Schur(\Lambda)}\(E(\lambda),E(\beta)\)$ to
$\End_{\H}\(N(\lambda),N(\beta)\)$; by Proposition~\ref{S' basis}(i) and
Proposition~\ref{endo} this is an isomorphism. Consequently, if an indecomposable
tilting module $T(\lambda')$ is a direct summand of $E(\lambda)$ then
$\Sfun\(T(\lambda')\)$ is an indecomposable direct summand of
$N(\lambda)$. Therefore,  by Proposition~\ref{twisted Youngs},
$\Sfun\(T(\lambda')\)\cong Y'(\mu)$ for some multipartition $\mu$.
Now,
$E(\lambda)\cong T(\lambda')\oplus
            \bigoplus_{\lambda'\gdom\mu}T(\mu)^{e_{\lambda\mu}}$
by Theorem~\ref{tilting modules}(i) and 
$$\Sfun\(E(\lambda)\)\cong N(\lambda)
   \cong Y'(\lambda)\oplus_{\nu\gdom\lambda}Y'(\nu)^{c_{\lambda\mu}}$$
by Proposition~\ref{twisted Youngs}(i). Hence, the result follows by induction on
the dominance ordering.
\end{proof}

\begin{Corollary}
Suppose that $R$ is a field,  $\omega\in\Lambda$
and that $Q_1,\dots,Q_r$ are distinct non--zero elements of $R$. Let
$\lambda$ and~$\mu$ be multipartitions of~$n$. Then
$[T(\lambda'):\nabla(\mu')]=[\W(\mu):L(\lambda)]$.
\label{tilting mults}
\end{Corollary}

\begin{proof}
Now $\Sfun\(\W(\nu)\)\cong S(\nu)$ by (\ref{Schur functor})(i); 
therefore, as $\Sfun$ projects onto the $\omega$--weight space,
$\Sfun\(\nabla(\nu)\)=\Sfun\(\W(\nu)^\dual\)\cong S(\nu)^\dual$.
Consequently, we have
\begin{xalignat*}{3}
[T(\lambda'):\nabla(\mu')]
  &=[\Sfun\(T(\lambda')\):\Sfun\(\nabla(\mu')\)],\\
  &=[Y'(\lambda):S(\mu')^\dual]
&&\text{by Proposition~\ref{tilting Youngs}},\\
  &=[Y'(\lambda):S'(\mu)],
&&\text{by Corollary~\ref{dual Spechts}},\\
  &=[\W'(\mu):L'(\lambda)],
          &&\text{by Proposition~\ref{twisted Youngs}(iv)},\\
  &=[\W(\mu):L(\lambda)],
\end{xalignat*}
where the last equality follows because the isomorphism
$\Schur(\Lambda)\cong\Schur'(\Lambda)$ of
Proposition~\ref{S' basis}(iii) identifies the Weyl modules $\W(\mu)$ and $\W'(\mu)$,
and the simple modules $L(\lambda)$ and $L'(\lambda)$.
\end{proof}

To conclude, we remark that in the case of the $q$--Schur algebras
(i.e.~when $r=1$) our proof of Corollary~\ref{tilting mults} looks quite
different to Donkin's~\cite{Donkin:tilt}; however, in spirit the two
arguments are the same in that they both rely on a duality between the
symmetric and exterior powers and on the isomorphism of 
Proposition~\ref{S' basis}(iii).

\section*{Acknowledgements}

I would like to thank Steve Donkin for explaining the 
tilting modules of the $q$--Schur algebras to me and for answering my 
numerous questions about them.

\let\em\it


\end{document}